\documentclass{amsart} 
\usepackage[mathscr]{eucal}
\usepackage{amssymb, amsmath,array, amscd}
\usepackage[OT2,T1]{fontenc}

\newtheorem{thm}{Theorem}

\newtheorem{prop}{Proposition}

\newtheorem{lemma}{Lemma}[section]
\newtheorem{claim}{Claim}[section]

\newtheorem{defi}{Definition}

\newtheorem{rem}{Remark}[section]

\newtheorem{cor}{Corollary}[section]

\newtheorem{Cor}{Corollary}

\newtheorem{prob}{Problem}

\newtheorem{ex}{Example}

\def\wt{\widetilde}
\def\dd{\mathcal{D}}
\def\wh{\widehat}

\def\fol{\mathbb{F}\text{ol}}

\def\p{\mathbb{P}}

\def\Te{\Theta}
\def\te{\theta}
\def\lim{\underset{z\to 0}{lim}\,}
\def\pa{\partial}
\def\sup{\supset}
\def\sub{\subset}
\def\De{\Delta}
\def\I{\mathcal{I}}
\def\d{\delta}

\def\emp{\emptyset}
\def\ov{\overline}
\def\om{\omega}
\def\Om{\Omega}
\def\fa{\mathcal{F}}
\def\a{\alpha}
\def\be{\beta}
\def\R{\mathbb{R}}
\def\ep{\epsilon}
\def\C{\mathbb{C}}

\def\{{\lbrace}
\def\}{\rbrace}

\def\la{\lambda}

\def\H{\mathcal{H}}

\def\g{\gamma}

\def\N{\mathbb{N}}
\def\Si{\Sigma}
\def\te{\theta}
\def\si{\sigma}
\def\G{\mathcal G}

\def\Z{\mathbb{Z}}
\def\O{\mathcal{O}}
\def\*{\star}

\def\Q{\mathbb{Q}}

\def\var{\varphi}

\def\X{\mathcal{X}}

\begin{document}

\title{The isotropy group of a foliation: the local case.}

\author{D. Cerveau \& A.  Lins Neto}

\begin{abstract}
Given a holomorphic singular foliation $\fa$ of $(\C^n,0)$ we define $Iso(\fa)$ as the group of germs of biholomorphisms on $(\C^n,0)$ preserving $\fa$: $Iso(\fa)=\{\Phi\in Diff(\C^n,0)\,|\,\Phi^*(\fa)=\fa\}$. The normal subgroup of $Iso(\fa)$, of biholomorphisms sending each leaf of $\fa$ into itself, will be denoted as $Fix(\fa)$.
The corresponding groups of formal biholomorphisms will be denoted as $\wh{Iso}(\fa)$ and $\wh{Fix}(\fa)$, respectively.  
The purpose of this paper will be to study the quotients $Iso(\fa)/Fix(\fa)$ and
$\wh{Fix}(\fa)/\wh{Fix}(\fa)$, mainly in the case of codimension one foliation.
\end{abstract}



\maketitle

\tableofcontents

\section{Introduction}\label{ss:1}
Let $\fa$ be a germ at $0\in\C^n$ of a singular codimension $p$ holomomorphic foliation, where $1\le p\le n-1$.
It is known that $\fa$ can be defined by a germ of holomorphic p-form $\Om\in\Om^p(\C^n,0)$, which is integrable in the following sense:
\begin{itemize}
\item[(1).] $\Om$ is locally completely decomposable outside its singular set: if $m\notin Sing(\Om)$, that is $\Om(m)\ne0$, then the germ $\Om_m$, of $\Om$ at $m$, is completely decomposable; i. e. there exist germs of holomorphic 1-forms at $m$, say $\om_1,...,\om_p$, such that $\Om_m=\om_1\wedge...\wedge\om_p$.
\item[(2).] $\om_1,...,\om_p$ satisfy the Frobenius integrability condition: $d\om_j\wedge\Om\equiv0$, $\forall$ $1\le j\le p$.   
\end{itemize}

Condition (1) implies that we can define a codimension $p$ distribution $\dd$ outside $Sing(\Om)$ by
\[
\dd(m)=\{v\in T_m\R^n\,|\,i_v\Om(m)=0\}:=ker(\Om(m))\,.
\]
Condition (2) implies that the distribution $\dd$ is integrable.
\vskip.1in
\begin{rem}
{\rm If $Cod_\C(Sing(\Om))\ge2$ and $\wt\Om$ is another germ of $p$-form that represents $\fa$ then there exists an unity $u\in\O_n^*$ such that $\wt\Om=u.\,\Om$.
From now on, we will assume that $Sing(\fa)=Sing(\Om)$ has codimension $\ge2$.}
\end{rem}
Associated to the foliation $\fa$ we introduce the following group
\begin{equation}\label{eq:1}
Iso(\fa)=\{\phi\in Diff(\C^n,0)\,|\,\phi^*(\fa)=\fa\}\,,
\end{equation}
the subgroup of $Diff(\C^n,0)$ of germs of holomorphic diffeomorphisms fixing the foliation $\fa$.

In terms of a germ of $p$-form $\Om$ defining $\fa$ relation (\ref{eq:1}) means that $\phi^*(\Om)$ also defines the foliation $\fa$, and so $\phi^*(\Om)=u.\,\Om$, where $u\in\O_n^*$.

We will consider also
\[
\wh{Iso}(\fa)=\{\phi\in \wh{Diff}(\C^n,0)\,|\,\phi^*(\Om)=\wh{u}.\,\Om\,,\,\wh{u}\in\wh\O_2^*\}\,,
\]
the group of formal diffeomorphisms "fixing" $\fa$. 

\begin{rem}
{\rm Let $\fa$ and $\Om$ be as before. The reader can check that a germ $\phi\in Iso(\fa)$ if, and only if, there are neighborhoods $U$ and $V$ of $0\in\C^n$, and representatives of $\Om$ and $\phi$, denoted by the same symbols, such that $\Om\in\Om^p(U\cup V)$, $\phi\colon U\to V$, and
\begin{itemize}
\item[a.] $\phi(Sing(\Om|_U))=Sing(\Om|_V)$.
\item[b.] $\phi$ sends leaves of $\fa|_U$ onto leaves of $\fa|_V$. 
\end{itemize}}
\end{rem}
Let us see some examples.
\begin{ex}\label{ex:1}
{\rm The simplest example of a germ $\fa$ of codimension $p$ foliation is the {\it regular foliation}, when $\Om(0)\ne0$. In this case, the definition and Frobenius theorem implies that there are local coordinates $z=(z_1,...,z_n)$ around $0\in\C^n$ and an unity $v\in\O_n^*$ such that $\Om=v(z)\,dz_1\wedge...\wedge dz_p$.
In this case, it is easy to see that $\phi=(\phi_1,...,\phi_n)\in Iso(\fa)$ if, and only if $\frac{\pa \phi_i}{\pa z_j}=0$, $\forall$ $1\le i\le p$, $\forall$ $j>p$.
In other words, if we set $\Psi_1=(\phi_1,...,\phi_p)$, $\Psi_2=(\phi_{p+1},...,\phi_n)$, $\zeta_1=(z_1,...,z_p)$ and $\zeta_2=(z_{p+1},...,z_n)$ then
\[
\phi(z)=(\Psi_1(\zeta_1),\Psi_2(\zeta_1,\zeta_2))\,.
\]}
\end{ex}
\begin{ex}\label{ex:2}
\rm Let $\fa$ and $\Om$ be as above. We say that a germ of holomorphic (resp. formal) vector field $X\in\X(\C^n,0)$ (resp.  $X\in\wh\X(\C^n,0)$) is tangent to the foliation $\fa$ if $i_X\Om=0$.
Let us denote the local flow of $X$ by $t\in(\C,0)\to X_t$. If $X$ is holomorphic and $t\in(\C,0)$ then $X_t$ defines a germ of diffeomorphism $X_t\colon(\C^n,0)\to(\C^n,X_t(0))$.
In particular, if $X(0)=0$ then $X_t\in Diff(\C^n,0)$, $\forall$ $t\in(\C,0)$.
The reader can verify directly that in this case $X_t\in Iso(\fa)$, $\forall$ $t\in(\C,0)$.

In this example, $X_t$ {\it fixes} the leaves of $\fa$, in the sense that if $L$ is a leaf of $\fa$ then $X_t(L)\cap L\ne\emp$ is an open subset of $L$.
\end{ex}

Motivated by example \ref{ex:2} we define $Fix(\fa)$ as the subgroup of $Iso(\fa)$ of germs of diffeomorphisms that "fix the leaves of $\fa$".
In other words, if $\Phi\in Fix(\fa)$, $U$ is the domain of $\Phi$ and $L\sub U$ is a leaf of $\fa|_U$ then $\Phi(L)\cap L\ne\emp$ is an open subset of $L$. For instance, if $L$ accumulates in the origin $0\in\C^n$ then, as a germ, we have $\Phi(L)=L$.

The formal complection of $Fix(\fa)$ in $\wh{Diff}(\C^n,0)$ will be denoted by $\wh{Fix}(\fa)$.
\begin{rem}
{\rm $Fix(\fa)$ (resp. $\wh{Fix}(\fa)$) is a normal subgroup of $Iso(\fa)$ (resp. $\wh{Iso}(\fa)$). We leave the proof to the reader.}
\end{rem}

One of the goals of this paper is to describe the groups $Iso(\fa)$ (resp. $\wh{Iso}(\fa)$) and the quotient $Iso(\fa)/Fix(\fa)$ (resp. $\wh{Iso}(\fa)/\wh{Fix}(\fa)$) in certain cases.
Let us begin by a simple example.

\begin{ex}
\rm Let $\om\in \Om^1(\C^n,0)$ be an integrable 1-form and
$\Phi\in Diff(\C^n,0)$ be an involution ($\Phi^2=I$), which is not in $Iso(\fa_\om)$.
In particular, $\eta:=\om\wedge\Phi^*(\om)\ne0$ and if $n\ge3$ then $\eta$ is integrable and defines a foliation $\fa_\eta$ of codimension two of $(\C^n,0)$. Since $\Phi^2=I$ we have
\[
\Phi^*(\eta)=\Phi^*(\om\wedge\Phi^*(\om))=\Phi^*(\om)\wedge\om=-\eta\,\,\implies\,\,
\Phi\in Iso(\fa_\eta)\,.
\] 
We would like to observe that $\Phi\notin Fix(\fa_\eta)$, so that its class in
$Iso(\fa_\eta)/Fix(\fa_\eta)$ is non-trivial.
\end{ex}

Another simple example is the following:

\begin{ex}\label{ex:hom}
{\rm We say that $\fa$ is a {\it homogeneous}\, foliation on $\C^n$ if there exists a p-form $\Om$ defining $\fa$ with all coefficients homogeneous of the same degree.
Let $R=\sum_jz_j\frac{\pa}{\pa z_j}$ be the radial vector field. The $p$-form $\Om$ is homogeneous if, and only if, $L_R\Om=(d+p)\,\Om$, where $d$ is the degree of the coefficients.
In particular, the group of dilatations of $\C^n$
\[
H=\{h_\rho(z)=\rho.\,z\,|\,\rho\in\C^*\}
\]
is contained in $Iso(\fa)$: $h_\rho^*(\Om)=\rho^{d+p}.\,\Om$ and so $h_\rho\in Iso(\fa)$.

We note however that $i_R\Om\equiv 0$ if, and only if, $H\sub Fix(\fa)$.
When $i_R\Om=0$ we say that the foliation $\fa_\Om$ is {\it conic}.

If $i_R\Om\not\equiv0$ then the generic element of $H$ does not fix the leaves of $\fa$ and the quotient $H/H\cap Fix(\fa)$ can be very complicated in general, as in the next example.}
\end{ex}

\begin{ex}\label{ex:5}
{\rm Let $\fa$ be the germ of codimension one logarithmic foliation given by $\Om=x_1...x_n\sum_{j=1}^n\la_j\frac{dx_j}{x_j}$, where $\la_1,...,\la_n\in\C$ are linearly independent over $\Z$. In this case the holonomy group of the leaf $L:=(x_1=0)\setminus\bigcup_{j>1}(x_j=0)$ in the transversal section $\Si:=\{(\tau,1,...,1)\,|\,\tau\in\C\}$ is linear and its multipliers are in the sub-group $G$ of $\C^*$ generated by $e^{2\pi i\la_j/\la_1}$, $2\le j\le n$. It can be shown that $H/H\cap Fix(\fa)$ is isomorphic to $\C^*/G$.}
\end{ex}

\begin{ex}\label{ex:6}
\rm Let $\om_1,...,\om_p$ be $p$ germs of integrable 1-forms on $(\C^n,0)$, where $2\le p<n$. 
Let $\eta:=\om_1\wedge...\wedge\om_p$. We will assume that $\eta\not\equiv0$. Note that $\eta$ is integrable.
We will assume also that:
\begin{itemize}
\item[(1).] If $\om=\sum_{j=1}^pf_j.\,\om_j$ is integrable, where $f_1,...,f_n\in\O_n$, and there is i such that $f_i\in\O_n^*$ then $\om=f_i.\,\om_i$.
\item[(2).] If $i\ne j$ then there is no $\Phi\in Iso(\fa_\eta)$ such that $\om_i\wedge\Phi^*(\om_j)=0$.
\end{itemize}
With the above hypothesis we have $Iso(\fa_\eta)=\bigcap_{j=1}^pIso(\fa_{\om_j})$. 

In fact, if $\Phi\in Iso(\fa_\eta)$ then $\Phi^*(\eta)=u.\,\eta$, $u\in\O^*_n$, implies that
$\Phi^*(\om_j)=\sum_{i=1}^pf_{ji}.\,\om_i$, $\forall1\le j\le p$, where $f_{ji}\in\O^*_n$ for some
$i$. Therefore, from (1) we get $\Phi^*(\om_j)=f_{ji}.\,\om_i$ and by (2) we get $i=j$ and
$\Phi^*(\om_j)=f_{jj}.\,\om_j$, so that $\Phi\in Iso(\fa_{\om_j})$, $\forall j$. The converse statement is immediate and is left to the reader.

A concrete example of this situation is the following: let $f,g\in\O_2$ be such that $f^{-1}(0)$ and $g^{-1}(0)$ are not homeomorphic. Define $\om_1,\om_2\in\Om^1(\C^4,0)$ as
$\om_1=z_1\,dz_2+df(z_1,z_2)$ and $\om_2=z_3\,dz_4+dg(z_3,z_4)$.

Using that $d\om_1=dz_1\wedge dz_2$ and $d\om_2=d_3\wedge dz_4$ it is possible to prove that $\om_1$ and $\om_2$ satisfy hypothesis (1) and (2) above. In particular, we get $Iso(\fa_{\om_1\wedge\om_2})=$
$=\{(\phi(z_1,z_2),\psi(z_3,z_4))\,|\,\phi\in Iso(\fa_{\om_1})$ and $Iso(\fa_{\om_2})\,\}\simeq Iso(\fa_{\om_1})\times Iso(\fa_{\om_2})$.
\vskip.1in

Another example in $(\C^3,0)$, satisfying (1) and (2), is the following: let $\om_1=a(x,y)\,dx+b(x,y)\,dy$ and $\om_2=dz$, where we assume that $\om_1$ has no meromorphic integrating factor. Here we have $Iso(\fa_{\om_1\wedge\om_2})=$
\[
=\{(\phi(x,y),\psi(z))\,|\,\phi\in Iso(\fa_{\om_1})\,\,\text{and}\,\,\psi\in Diff(\C,0)\}\simeq Iso(\fa_{\om_1})\times Diff(\C,0)\,.
\]
\end{ex}

We begin in \S\,\ref{ss:2} studying the case of codimension one foliations. In \cite{BCM} the authors study some special cases of this situation: 
when the dimension is $n=2$ and when $\fa$ has an holomorphic first integral, related with a known Briancon-Skoda theorem \cite{BCM}.
In particular, we intend here to precise and generalize some of the results of \cite{BCM}.
Our first result in this direction is the following:
\begin{thm}\label{t:1}
Let $\fa$ be a germ at $0\in\C^n$ of codimension one foliation defined by a germ of integrable 1-form $\Om$. Suppose that the quotient $\wh{Iso}(\fa)/\wh{Fix}(\fa)$ has an element of infinite order.
Then $\fa$ is formally Liouville integrable, that is, there exists $\hat{f}\in\wh\O_n$ such that $\frac{1}{\wh{f}}\,\Om$ is closed.
\end{thm}

In \S\, \ref{ss:23} we will apply theorem \ref{t:1} in the case of foliations on $(\C^2,0)$ that in the process of resolution have a non-dicritical irreducible component with non-abelian holonomy. In theorem \ref{t:2} we will prove that for such a foliation $\fa$ then $\wh{Iso}(\fa)/\wh{Fix}(\fa)$ is isomorphic to a finite subgroup of the linear group $GL(2,\C)$.
\vskip.1in
Theorem \ref{t:1} is stated in a formal context. A germ of meromorphic closed 1-form ${\om}/{f}$, $\om\in\Om^1(\C^n,0)$, $f\in\O_n$, can be written in a normal form of the following type: if $f=f_1^{n_1+1}...f_r^{n_r+1}$ is the decomposition of $f$ into irreducible factors, then it can be proved that (cf. \cite{CeMa})
\begin{equation}\label{eq:fe}
\frac{\om}{f}=\sum_{j=1}^r\la_j\frac{df_j}{f_j}+d\left(\frac{H}{f_1^{n_1}...f_r^{n_r}}\right)\,,
\end{equation}
where $H\in\O_n$ and $\la_j\in\C$, $1\le j\le r$.

In the formal case, when $\om\in\wh{\Om}^1(\C^n,0)$ and $f_j\in \wh\O_n$, $1\le j\le r$, there is a similar result, proved in \cite{th}, with $H\in\wh\O_n$ in (\ref{eq:fe}).
In our case, $\om\in\Om^1(\C^n,0)$ and $f\in\wh\O_n$, but there are conditions, on the residues $\la_j$, assuring the convergence of $f$ (see \cite{CeMa}).

About the cardinality of the quotient group and the existence of holomorphic integrating factor, by using results of P\'erez Marco \cite{pm} it is possible to construct examples of germs of foliations $\fa$ on $(\C^2,0)$ in which $Iso(\fa)/Fix(\fa)$ is non-countable, contains elements of infinite order, but $\fa$ is not liouville integrable: it cannot be defined by a meromorphic closed form (see theorem \ref{t:3} in \S\,\ref{ss:24}).

In \S\,\ref{ss:3} we will study $Iso(\fa)$ and $Fix(\fa)$ when $\fa$ is homogeneous or has an integrating factor.

\begin{defi}
\rm We say that a p-form $\Om$ on $\C^n$ is conical if it is homogeneous and $i_R\Om=0$, where $R$ is the radial vector field on $\C^n$. A holomorphic foliation of codimension $p$ on $\C^n$ is conic if it can be defined by a conical p-form.
\end{defi}

If $\fa$ is a conical foliation on $\C^n$ , as above, then $\fa$ induces a foliation $\wt\fa$, of the same codimension, on the projective space $\p^{n-1}$. In theorem \ref{t:4} we prove that, in this case, $Iso(\fa)/Fix(\fa)$ is isomorphic to $Aut(\wt\fa)$, the subgroup of $Aut(\p^{n-1})$ of automorphisms of $\p^{n-1}$ preserving $\wt\fa$. Let us see an example.
\begin{ex}\label{ex:5}
{\rm Jouanolou's example of degree $d\ge2$ on $\p^n$ (see \cite{j} and \cite{LS}) is defined in homogeneous coordinates of $\C^{n+1}$ by the homogeneous $(n-1)$-form $\Om=i_Ri_X\nu$, where $R$ is the radial vector field on $\C^{n+1}$, $\nu=dx_1\wedge...\wedge dx_{n+1}$ and
\[
X=x_{n+1}^d\frac{\pa}{\pa x_{1}}+\sum_{j=2}^nx_{j-1}^d\frac{\pa}{\pa x_j}\,.
\]
The foliation $\wt\fa$, induced on $\p^n$ by $\fa_\Om$, can be defined in the affine coordinate $x_{n+1}=1$ by the vector field
\[
Y=(1-x_n^d.\,x_1)\frac{\pa}{\pa x_1}+\sum_{j=2}^{n}(x_{j-1}^d-x_n^d.\,x_j)\frac{\pa}{\pa x_j}\,.
\]

Let $D=\frac{d^{n+1}-1}{d-1}$ and $\la$ be a primitive $D^{th}$-root of unity.
It is known that $Aut(\wt\fa)$ is isomorphic to the finite sub-group $G(n,d)$ of $PSL(n+1,\C)$ generated by the transformations $\tau$ and $\rho$, where
\[
\tau(x_1,...,x_n,x_{n+1})=\left(\la.\,x_1,\la^{d+1}.\,x_2,...,\la^{\frac{d^j-1}{d-1}}.\,x_j,...,\la^{\frac{d^n-1}{d-1}}.\,x_n,x_{n+1}\right)
\]
and
\[
\rho(x_1,x_2,...,x_n,x_{n+1})=(x_{n+1},x_1,x_2,...,x_n)\,.
\]
In particular, $Iso(\fa_\Om)/Fix(\fa_\Om)\simeq G(n,d)$ by theorem \ref{t:4}.}
\end{ex}
\vskip.1in

On the other hand, if $\om\in\Om^1(\C^n)$ is integrable and homogeneous, but is non-conical, then the following facts are known (see \cite{CeMa}):
\begin{itemize}
\item[a.] If $i_R\om =f\ne0$ then $f$ is an integrating factor of $\om$: $d\left(f^{-1}.\,\om\right)=0$.
\item[b.] $f$ is homogeneous: $R(f)=\ell\,f$, for some $\ell\in\N$.
\item[c.] If $f=f_1^{n_1+1}...f_r^{n_r+1}$ is the decomposition of $f$ into homogeneous factors then the hypersurfaces $(f_j=0)$ are $\fa_\om$-invariant, $\forall$ $1\le j\le r$. Moreover, the form $f^{-1}.\,\om$ can be written as in (\ref{eq:fe}).
\item[d.] When $f$ is reduced, $f=f_1...f_r$, then $\Om:=f^{-1}.\,\om$ is logarithmic:
\begin{equation}\label{eq:log}
\Om=\sum_{j=1}^r\la_j\frac{df_j}{f_j}\,.
\end{equation}
\end{itemize}

When $\Om=f^{-1}.\,\om$ is like in (\ref{eq:fe}) or (\ref{eq:log}) and $f_1,...,f_r\in\O_n$ are relatively prime, but not necessarily homogeneous, we will consider in \S\,\ref{ss:33} the case in which $\fa_\Om$ has no meromorphic first integral. In proposition \ref{p:4} we will prove that if $\Om$ is a closed meromorphic 1-form, without meromorphic first integral, and if $\Phi\in Iso(\fa_\Om)$ then $\Phi^*(\Om)=\d\,\Om$, where $\d$ is a root of unity.
As we will see there, in this case $\Phi$ permutes the "separatrices" $(f_j=0)$, and if $\Phi(f_j=0)=(f_j=0)$, $1\le j\le r$, then $\d=1$.

In theorem \ref{t:5} we will prove that if $\Om$ is logarithmic, has no meromorphic first integral and $Iso(\fa_\Om)$ contains a transformation $\Phi$ with $D\Phi(0)=\rho.\,I$, where $\rho$ is not a root of unity, then there exists $h\in Diff(\C^n,0)$ such that $h^*\,\Om$ is homogeneous.

In corollaries \ref{c:2}, \ref{c:3}, \ref{c:4} and \ref{c:5}, we apply theorem \ref{t:5} in some special situations.

\begin{ex}\label{ex:c4}
\rm Let $f_1,f_2\in \O_2$ be such that:
\begin{itemize}
\item[1.] $f_1$ is irreducible in $\O_n$ and its first non-vanishing jet at $0\in\C$ is $J^3_0(f_1)=z_1.\,z_2^2$.
\item[2.] $J^1_0(f_2)=z_1+z_2$.
\end{itemize}
Let $\Om=\frac{df_1}{f_1}+\la\,\frac{df_2}{f_2}$, where $\la\notin\Q$. Then $Iso(\fa_\Om)=Fix(\fa_\Om)$.

This example is an application of corollary \ref{c:4}, as we will see in \S\,\ref{ss:31}.
\end{ex}

\begin{rem}
\rm Let $\om\in \Om^1(\C^n,0)$ be an integrable form inducing the germ of foliation $\fa_\om$. Note that $Iso(\fa_\om)/Fix(\fa_\om)$ can be considered as a subgroup of $\wh{Iso}(\fa_\om)/\wh{Fix}(\fa_\om)$. 

In fact, there is a natural group inclusion $In\colon Iso(\fa_\om)/Fix(\fa_\om)\to \wh{Iso}(\fa_\om)/\wh{Fix}(\fa_\om)$, as the reader can check. This map is injective, but not surjective in general.

For instance, in the case of 1-forms $\om\in \Om^1(\C^2,0)$ that are formally linearizable, but not holomorphically linearizable, studied in theorem \ref{t:3}, the map is not surjective.
In contrast, in the case of conical foliations, to be studied in \S\,\ref{ss:32}, the map is an isomorphism.
A natural problem is the following:
\begin{prob}
\rm When $In\colon Iso(\fa_\om)/Fix(\fa_\om)\to \wh{Iso}(\fa_\om)/\wh{Fix}(\textsl{}\fa_\om)$ is an isomorphism?
\end{prob}
\end{rem}

Another natural problem is the following:

\begin{prob}\label{pr:2}
\rm Let $\fa$ be a germ of foliation such that any element $\Phi\in {Iso}(\fa)/{Fix}(\fa)$ has finite order.
Is ${Iso}(\fa)/{Fix}(\fa)$ finite? The same question can be posed for
$\wh{Iso}(\fa)/\wh{Fix}(\fa)$.
\end{prob}

In theorem \ref{t:2} we prove that the answer of problem \ref{pr:2} is positive in a particular case in dimension two. In corollary \ref{c:4} of theorem \ref{t:5} we prove a similar statement in the case of logarithmic foliations in $(\C^n,0)$ (see \S\,\ref{ss:31}). 

\begin{rem}\label{r:15}
\rm Let $S=\sum_{j=1}^nk_j\,z_j\frac{\pa}{\pa z_j}$ be a semi-simple linear vector field with eigenvalues $k_1,...,k_n\in\Z_{>0}$ and $gcd\,(k_1,...,k_n)=1$. We say that a $p$-form
$\eta\in\Om^p(\C^n,0)$ is $S$ quasi-homogeneous if $L_S\eta=k\,\eta$, where $k\in\N$. 
In the case of a germ of function $f\in\O_n$ the identity $S(f)=L_S(f)=k.\,f$ means that
\[
f(\tau^{k_1}z_1,...,\tau^{k_n}z_n)=\tau^kf(z_1,...,z_n)\,\,,\,\,\forall (z_1,...,z_n)\in(\C,0)\,\,,\,\,\forall \tau\in\C\,.
\]

We would like to observe the following facts:
\begin{itemize}
\item[1.] $\eta$ is homogeneous if, and only if, $\eta$ is quasi-homogeneous with respect to $R$, the radial vector field in $\C^n$.
\item[2.] If $\eta$ is $S$ quasi-homogeneous then their coefficients are $S$ quasi-homogeneous. As a consequence all coefficients of $\eta$ are polynomials.
\item[3.] The flow $exp(t\,S)$ of $S$ induces a $\C^*$ action $\Psi\colon \C^*\times \C^n\to\C^n$:
\[
\Psi_\tau(z_1,...,z_n):=\Psi(\tau,z_1,...,z_n)=(\tau^{k_1}z_1,...,\tau^{k_n}z_n)\,.
\]
\end{itemize}
The relation $L_S\eta=k\eta$ is equivalent to $\Psi_\tau^*(\eta)=\tau^k\eta$.

We say that $\eta$ is a $S$ conical foliation if $\eta$ is $S$ quasi-homogeneous and $i_S\eta=0$.
A $S$ conical and integrable form induces a foliation on the weighted projective space $\p^{n-1}_k$, associated to the weights $k=(k_1,...,k_n)$ (cf. \cite{dg}).

A natural question is the following:

\begin{prob}
\rm Let $S=\sum_{j=1}^nk_j\,z_j\frac{\pa}{\pa z_j}$ be as above. Are there statements concerning $S$ conical and $S$ quasi-homogeneous foliations, similar to theorems \ref{t:4} and \ref{t:5}?
\end{prob}
\end{rem}

\section{Theorem \ref{t:1} and correlated facts}\label{ss:2}

\subsection{Preliminaries}\label{ss:21}
This section is devoted to the statement of some well known results that will be used along the text.
We will prove also theorem \ref{t:1} in \S\,\ref{ss:22} and give an application in \S\,\ref{ss:23}.
In \S\,\ref{ss:24} we construct examples of foliations $\fa$ in dimension two for which $Iso(\fa)/Fix(\fa)$ is non-countable. These examples are not Liouville integrable, but formally Liouville integrable, in the sense that they admit a formal integrating factor.

We begin recalling the concept of unipotent germ of formal diffeomorphism.
\begin{defi}
\rm We say that $\phi\in \wh{Diff}(\C^n,0)$ is unipotent if the linear part $D\phi(0)$ is unipotent.
A formal vector field $X\in\wh\X(\C^n,0)$ is nilpotent if its linear part $DX(0)$ is nilpotent.
\end{defi}
In fact, the following proposition is known (cf. \cite{AR} and \cite{BCM}):
\begin{prop}\label{p:21}
A formal germ $\phi\in\wh{Diff}(\C^n,0)$ is unipotent if, and only if, there exists a nilpotent formal vector field $X\in\wh\X(\C^n,0)$ such that $\phi=exp\,(X)$, where
$t\mapsto exp(t\,X)$ denotes the formal flow of $X$.
Moreover, if $X\in\wh\X(\C^n,0)$ is nilpotent then the formal flow is polynomial in $t$, in the sense that
\begin{equation}\label{eq:2}
exp\,(t\,X)(z)=\sum_\si\,P_\si(t)\,z^\si\,,
\end{equation}
where $P_\si(t)$ is a polynomial in $t$ for all $\si=(\si_1,...,\si_n)$.
\end{prop}
As a consequence, we have the following:
\begin{cor}\label{c:21}
If $\phi\in\wh{Iso}(\fa)$ is unipotent and $\phi=exp\,(X)$ then $\phi_t=exp\,(t\,X)\in\wh{Iso}(\fa)$, $\forall t\in\C$. 
\end{cor}
{\it Proof.}
Let $\Om$ be a germ of integrable $p$-form that represents $\fa$. Since $\phi\in \wh{Iso}(\fa)$ we have
\[
\phi^*(\Om)=\wh{u}.\,\Om\,,\,\wh{u}\in\O_n^*\,.
\]
From the above relation, we get by induction on $k\in\N$ that
\begin{equation}\label{eq:3}
(\phi^k)^*(\Om)=\wh{u}_k.\,\Om\,,
\end{equation}
where $\wh{u}_{k}$ is defined inductively by $\wh{u}_1=\wh{u}$ and $\wh{u}_{k+1}=\wh{u}_k\circ\phi.\,\wh{u}$, if $k\ge1$.

Write
\[
\Om=\sum_IF_I(z)\,dz^I\,,
\]
where $I=(1\le i_1<...<i_p\le n)$, $dz^I=dz_{i_1}\wedge...\wedge dz_{i_p}$ and $F_I(z)=\sum_\si A_{I,\si}z^\si$.

Let $X\in\wh\X(\C^n,0)$ be such that $\phi=exp(X)$.
Since $exp(t\,X)(z)=\sum_\si P_\si(t)\,z^\si$, where $P_\si\in\C[t]$, we get by direct substitution that
\[
exp(t\,X)^*(\Om)=\sum_IF_I(t,z)\,dz^I\,,
\]
where $F_I(0,z)=F_I(z)$ and
\begin{equation}\label{eq:4}
F_I(t,z)=\sum_\si Q_{I,\si}(t)\,z^\si\,\,,\,\,Q_{I,\si}\in\C[t]\,,\,\forall I,\si
\end{equation}
Now, we can write $\phi^k=exp(k\,X)$ if $k\in \Z$, so that, from (\ref{eq:3}) and (\ref{eq:4}) we obtain that for all $I_1,\si_1$ and $I_2,\si_2$ such that $A_{I_1,\si_1}.A_{I_2,\si_2}\ne0$ we have
\[
\frac{Q_{I_1,\si_1}(k)}{A_{I_1,\si_1}}=\frac{Q_{I_2,\si_2}(k)}{A_{I_2,\si_2}}\,,\,\forall k\in \Z\,.
\]
Since all $Q_{I,\si}(t)$ are polynomials in $t$, from the above relation that
\[
\frac{Q_{I_1,\si_1}(t)}{A_{I_1,\si_1}}=\frac{Q_{I_2,\si_2}(t)}{A_{I_2,\si_2}}\,,\,\forall t\in \C
\]
which implies the corollary.
\qed
\vskip.1in
Another well known fact is the following:
\begin{prop}\label{p:22}
Let $\phi\in\wh{Diff}(\C^n,0)$ be a formal diffeomorphism. Then $\phi$ admits a formal Jordan decomposition: $\phi=\phi_S\circ\phi_U$, where $\phi_S$ and $\phi_U$ are commuting formal diffeomorphisms, $\phi_S$ is semi-simple (formally conjugated to its linear part $D\phi_S(0)$) and $\phi_U$ is unipotent. 
\end{prop}

\begin{rem}\label{r:21}
{\rm If $\phi\in \wh{Iso}(\fa)$ then $\phi_S,\phi_U\in \wh{Iso}(\fa)$. The proof can be found in \cite{BCM}.}
\end{rem}

\subsection{Proof of theorem \ref{t:1} and complements}\label{ss:22}
Let $\fa$ be a germ of codimension one foliation defined by an integrable germ $\Om\in \Om^1(\C^n,0)$.
We will assume that $\wh{Iso}(\fa)$ has an element of infinite order in the quotient $\wh{Iso}(\fa)/\wh{Fix}(\fa)$, say $\phi$.
The proof will be based in several remarks. The first is elementary:
\begin{prop}\label{p:23}
Let $\phi\in \wh{Iso}(\fa)$ be unipotent and $X\in \wh\X(\C^n,0)$ be such that $\phi=exp(X)$. Let $f:=i_X\Om\in\wh\O_n$. We have two possibilities:
\begin{itemize}
\item[(a).] $f\not\equiv0$ and $\om:=\frac{1}{f}\,\Om$ is closed.
\item[(b).] $f\equiv0$ and in this case $\phi\in\wh{Fix}(\fa)$. 
\end{itemize}
\end{prop}
{\it Proof.}
From corollary \ref{c:21} we know that $exp(t\,X)\in \wh{Iso}(\fa)$, which means that
\[
exp(t\,X)^*\,\Om=u_t\,\Om\,,
\]
where $u_t\in \wh\O_n^*$ $\forall t\in\C$. Using that $\frac{d}{dt}exp(t\,X)^*\Om|_{t=0}=L_X\Om$, where $L_X$ denotes the Lie derivative in the direction of $X$, and taking the derivative in both members of the above relation we get
\[
L_X\Om=v\,\Om\,,\,v=\frac{d}{dt}\,u_t|_{t=0}\,.
\]
Since $v\,\Om=L_X\Om=i_X\,d\Om+d\,i_X\Om=i_X\,d\Om+df$ we get
\[
0=L_X\Om\wedge \Om=i_X\,d\Om\wedge\Om+df\wedge\Om\,\,\implies\,\,df\wedge\Om=\Om\wedge i_X\,d\Om\,.
\]
On the other hand, the integrability condition $\Om\wedge d\Om=0$ implies that
\[
0=i_X(\Om\wedge d\Om)=i_X\Om.\,d\Om-\Om\wedge i_X\,d\Om=f\,d\Om-\Om\wedge i_X\,d\Om\,\,\implies
\]
\[
f\,d\Om=\Om\wedge i_X\,d\Om=df\wedge\Om\,\,\iff\,\,d\left(\frac{1}{f}\,\Om\right)=0\,\,\text{, if $f\ne0$}.
\]
This proves (a).

If $X$ is holomorphic and $i_X\Om=0$ then the orbits of the flow $exp(t\,X)$ are contained in the leaves of $\fa$, so that $\phi=exp(X)\in Fix(\fa)$.

In the general case the formal flow $t\mapsto exp(t\,X)$ is tangent to $\fa$.
\qed
\vskip.1in
Let us finish the proof of theorem \ref{t:1}. Let $\fa$ be defined by the germ of integrable 1-form $\Om$ and $\phi\in\wh{Iso}(\fa)$ be of infinite order in $\wh{Iso}(\fa)/\wh{Fix}(\fa)$.
By proposition \ref{p:22} we can decompose $\phi=\phi_S\circ\phi_U$, where $\phi_S$ is semi-simple and $\phi_U$ unipotent. By remark \ref{r:21} $\phi_S$ and $\phi_U$ are in $\wh{Iso}(\fa)$. Let $X\in\X(\C^n,0)$ be such that $\phi_U=exp(X)$. If $i_X\,\Om\not\equiv0$ we are done by proposition \ref{p:23}.
If $i_X\Om\equiv0$ then $\phi_U\in\wh{Fix}(\fa)$ and this implies that the classes of $\phi_S$ and $\phi$ in $\wh{Iso}(\fa)/\wh{Fix}(\fa)$ coincide. Hence, $\phi_S$ is of infinite order in $\wh{Iso}(\fa)/\wh{Fix}(\fa)$. Now, $\phi_S$ is linearizable: there exists $\var\in \wh{Diff}(\C^n,0)$ such that $\var^{-1}\circ \phi_S\circ \var=L$, where $L=D\phi_S(0)$ is linear and semi-simple: in some base of $\C^n$ the matrix of $L$ is diagonal and the subgroup of powers of $L$, $H=\{L^k\,|\,k\in\Z\}\sub GL(n,\C)$, is abelian.
Note that $L$ is not of finite order, for otherwise $\phi_S$ would be also of finite order. In particular $H$ is infinite.

Let $G$ be the Zariski closure of $H$. Since $H$ is infinite and abelian, $G$ is an abelian algebraic Lie group of dimension $\ge1$.
Let
\[
\wh{G}=\var\circ G\circ\var^{-1}=\{\var\circ g\circ\var^{-1}\,|\,g\in G\}\sub\wh{Diff}(\C^n,0)\,.
\]

Note that $\wh{G}$ is abelian.
We assert that $\wh{G}\sub\wh{Iso}(\fa)$.
In fact, let $\wt\Om=\var^*(\Om)$. From $\phi_S^*(\Om)=u.\,\Om$, where $u\in\wh\O_n^*$, we get 
$L^*(\wt{\Om})=\wt{u}.\,\wt\Om$ and $(L^k)^*(\wt\Om)=\wt{u}_k.\,\wt\Om$, $\wt{u}_k\in\wh\O_n^*$, $\forall k\in\Z$.
Since $G$ is the Zariski closure of $H$, for all $g\in G$ we must have $g^*(\wt{\Om})=\wt{u}_g.\,\wt\Om$, $\wt{u}_g\in\wh\O_n^*$.
Hence, if $\wh{g}=\var\circ g\circ\var^{-1}\in\wh{G}$ then $\wh{g}^*(\Om)=u_g.\,\Om$, $u_g=\wt{u}_g\circ\var^{-1}$, so that $\wh{g}\in \wh{Iso}(\fa)$.

Now, since $G$ is a complex Lie group of dimension $\ge1$ the Lie algebra $\G$ of $G$ has dimension $\ge1$. Moreover, if $X\in\G$ then $exp\,(t\,X)^*(\wt\Om)=u_t.\,\wt\Om$, so that
$L_X\wt\Om=h\,\wt\Om$, where $h=\frac{d}{dt}u_t|_{t=0}$.
We assert that there exists $X\in\G$ such that $i_X\wt\Om\not\equiv0$.

In fact, if not then $0\in\G$ has a neighborhood $V$ such that $i_X\wt\Om=0$ $\forall X\in V$. 
In particular, if $G^o$ is the connected component of $G$ containing the identity, then for any $g\in G^o$ we have $g=exp(t\,X)$ where $t\in\C$ and $X\in V$. This implies that the conjugate $\wh{g}:=\var\circ g\circ \var^{-1}\in \wh{Fix}(\fa)$.

On the other hand, since $G$ has a finite number of connected components there exists $k\in \N$ such that $L^k\in G^o$.
But, this would imply that the class of $\phi_S^k$ in $\wh{Iso}(\fa)/\wh{Fix}(\fa)$ would  be trivial, contradicting the hypothesis. 

Finally, if $X\in\G$ is such that $i_X\wt\Om\ne0$ then the vector field $Y=\var_*(X)\in\wh\X(\C^n,0)$ satisfies $f:=i_Y\Om\not\equiv0$, so that $d\left(\frac{1}{f}\Om\right)=0$.
\qed 

\subsection{An Application of theorem \ref{t:1}}\label{ss:23}
A consequence of theorem \ref{t:1} is that when a germ of codimension one foliation $\fa$, by some reason, cannot be represented by a germ of holomorphic form that has no formal integrating factor then all elements of $\wh{Iso}(\fa)/\wh{Fix}(\fa)$ are of finite order.
\vskip.1in
We will apply the above remark in the case of germs of foliations on $(\C^2,0)$.

Is is known that any germ at $0\in \C^2$ of foliation by curves, say $\fa$, has a resolution by a sequence of punctual blowing-ups (cf. \cite{se} and \cite{MM}).
After the resolution process $\Pi\colon(M,E)\to(\C^2,0)$, where $\Pi$ denotes the blowing-up map and $E\sub M$ the exceptional divisor, we obtain the resolved foliation $\Pi^*(\fa):=\wt\fa$. All the irreducible components of $E$ are biholomorphic to $\p^1$. Some of them, the non-dicritical components, are $\wt\fa$-invariant, whereas others, the dicritical ones, are not invariant. The foliation $\wt\fa$ has only simple singularities (see \cite{MM}). A singularity of a germ of holomorphic vector field $Z$ at $(\C^2,0)$ is simple if:
\begin{itemize}
\item[a.] The derivative $DZ(0)$ is semi-simple and not identically zero. Let $\la_1$ and $\la_2$ be the eigenvalues of $DZ(0)$.
\item[b.] If $\la_1\ne0$ and $\la_2=0$, or vice-versa, the singularity is called a saddle-node.
\item[c.] If $\la_1,\la_2\ne0$ then $\la_2/\la_1\notin\Q_+$. 
\end{itemize}

A non-dicritical component $D$ of $E$ contains necessarily singularities of $\wt\fa$, say $S=\{p_1,...,p_k\}$. If we fix a transverse section $\Si$ to $\wt\fa$ at a non-singular point $p\in D$ then the holonomy group of $\wt\fa$ at $D$ is a representation $\rho$ of $\Pi_1(D\setminus S)$ on the group of germs of biholomorphisms $Diff(\Si,p)\simeq Diff(\C,0)$.
The image $\rho(\Pi_1(D\setminus S))$ is called the holonomy group of $D$.

We would like to observe that, given a finitely generated subgroup $H$ of $Diff(\C,0)$ then it is possible to construct examples of foliations such that $\wt\fa$ has a non-dicritical divisor $D$ with holonomy group isomorphic to $H$ (see \cite{LN1}).
As an application of theorem \ref{t:1} we have the following result:

\begin{thm}\label{t:2}
Let $\fa$ be a germ at $0\in\C^2$.  Assume that after the resolution process, $\Pi\colon(M,E)\to(\C^2,0)$, the total divisor $E$ has a non-dicritical irreducible component $D$ for which the holonomy of the strict transform $\Pi^*(\fa)$ is non abelian.
Then $\wh{Iso}(\fa)/\wh{Fix}(\fa)$ is isomorphic to a finite subgroup of the linear group $GL(2,\C)$.
\end{thm}

{\it Proof.}
Assume that $\fa$ is defined by the germ of vector field $X=P\frac{\pa}{\pa x}+Q\frac{\pa}{\pa y}$, or equivalently, by the dual form $\Om=P\,dy-Q\,dx$, where $P,Q\in \O_2$ and $P(0)=Q(0)=0$.
\vskip.1in
\begin{lemma}\label{l:21}
$\Om$ has no formal integrating factor.
In particular, any $\phi\in\wh{Iso}(\fa)/\wh{Fix}(\fa)$ has finite order.
\end{lemma}

{\it Proof.}
Suppose by contradiction that $\Om$ has a formal integrating factor $\wh{f}$: $d\left(\frac{1}{\wh{f}}\,\Om\right)=0$. Let $\wh{f}=\Pi_{j=1}^kf_j^{r_j}$ be the decomposition of $\wh{f}$ into formal irreducible factors.
In this case, we can write (cf. \cite{CeMa} and \cite{th}):
\begin{equation}\label{eq:5}
\frac{1}{\wh{f}}\,\Om=\sum_{j=1}^k\la_j\frac{df_j}{f_j}+d\left(\frac{g}{f_1^{r_1-1}...f_k^{r_k-1}}\right)
\end{equation}
where $\la_j\in \C$ and $g\in\wh\O_2$.
Let $\Pi\colon (M,E)\to(\C^2,0)$ be the minimal resolution process of $\fa$ and $\fa^*$ be the strict transform of $\fa$ by $\Pi$.
Let $D\sub E$ be the irreducible component with non abelian holonomy group. 

\begin{claim}\label{cl:21}
\rm Let $p\in D\setminus Sing(\fa^*)$. Then there are formal coordinates $(t,\wh{x})$ around $p$ such that
\begin{itemize}
\item[a.] $(\wh{x}=0)\sub D$ and $p=(0,0)$.
\item[b.] $\Pi^*\left(\frac{1}{\wh{f}}\,\Om\right)=\phi(\wh{x})\,d\wh{x}$, where $\phi$ is of one of the following types:
\begin{itemize}
\item[1.] $\phi(\wh{x})=\wh{x}^m$, where $m\ge0$.
\item[2.] $\phi(\wh{x})=\la/\wh{x}$, where $\la\in\C^*$.
\item[3.] $\phi(\wh{x})=(1+\la\,\wh{x}^{\ell-1})/\wh{x}^\ell$, where $\la\in\C$ and $\ell\ge2$.
\end{itemize}
\end{itemize}
\end{claim} 

{\it Proof.}
Since $p\notin Sing(\fa^*)$ and $D$ is $\fa^*$-invariant there exists a holomorphic coordinate system $(U,(t,x))$ such that:
\begin{itemize}
\item[i.] $D\cap U=(x=0)$ and $p=(0,0)$.
\item[ii.] $\fa^*|_U$ is defined by $dx=0$, or equivalently their leaves are the levels $x=cte$.
\end{itemize}
Let $\wt{f}_j:=f_j\circ\Pi$ and $\wt{g}=g\circ\Pi$. Since $f_j(0)=0$ and $p\notin Sing(\fa^*)$ we must have $\wt{f}_j(t,x)=x^{m_j}.\,u_j(t,x)$, where $m_j\ge1$ and $u_j$ is a formal unity, $1\le j\le k$.
Since $\fa^*$ is defined by $dx=0$, we must have $\Pi^*(\Om)=v(t,x).\,x^\mu\,dx$, where $\mu\ge1$ and $v$ is a holomorphic unity.
In particular, we get
\[
\Pi^*\left(\frac{1}{\wh{f}}\,\Om\right)=\frac{\wh{v}}{x^\ell}.\,dx:=\wh{\phi}\,dx\,,
\]
where $\wh{v}=v/u_1^{r_1}...u_k^{r_k}$ is a formal unity, and $\ell=\sum_jr_j\,m_j-\mu\in\Z$.
Since $\frac{1}{\wh{f}}\,\Om$ is closed, we must have $d\wh{v}\wedge dx=0$, so that $\wh{v}=\wh{v}(x)$.

1. If $\ell=-m\le0$ then $\wh{\phi}(x)\,dx=\wh{v}(x).\,x^m\,dx$.
Let $H(x)=x^{m+1}.\,w(x)$ be a formal series with $H'(x)=x^m\,\wh{v}(x)$ and $w\in\wh\O_2^*$.
Since $w(0)\ne0$ there exists $u\in\wh\O_n^*$ with $u^{m+1}=(m+1)\,w$. If $\wh{x}=x\,u$ then $H=(m+1)^{-1}\,\wh{x}^{m+1}$ and $dH=\wh{x}^m\,d\wh{x}$.

2. If $\ell=1$ then $\wh{\phi}(x)\,dx$ has a pole of order one and we can write
\[
\wh{\phi}(x)=\frac{\la}{x}+\var(x)\,,
\]
where $\la\in\C^*$ and $\var(x)$ is a formal power series. If we set $u(x)=exp\,(\var(x)/\la)$ then $d\var=\la\frac{du}{u}$. In particular, if $\wh{x}=u(x)\,x$ then
\[
\wh{\phi}(x)\,dx=\la\frac{d\wh{x}}{\wh{x}}\,.
\]

3. If $\ell>1$ and $\wh{v}(x)=\sum_{j\ge1}a_j\,x^j$ then we can write
\[
\wh{\phi}(x)\,dx=\left(\frac{\la}{x}+\frac{\var(x)}{x^\ell}\right)\,dx=\frac{\var(x)+\la\,x^{\ell-1}}{x^\ell}\,\,dx\,,
\]
where $\la=a_{\ell-1}$ and $\var(0)=a_0\ne0$. In this case we consider the formal vector field $X=\frac{x^\ell}{\var(x)+\la\,x^{\ell-1}}\,\frac{\pa}{\pa x}$, for which $i_X(\wh{\phi}(x)\,dx)=1$. It is known that there exists $\wh{x}=\psi(x)\in\wh{Diff}(\C,0)$ such that $\psi^*(X)=\frac{\wh{x}^\ell}{1+\la\,\wh{x}^{\ell-1}}\,\frac{\pa}{\pa \wh{x}}$ (see \cite{ma} and \cite{CM}). It can be checked that $\psi^*(\wh{\phi}(x)\,dx)=\frac{1+\la\,\wh{x}^{\ell-1}}{\wh{x}^\ell}\,d\wh{x}$, which proves the claim.
\qed
\vskip.1in
\begin{claim}\label{cl:22}
\rm Let $h\in \wh{Diff}(\C,0)$ be such that $h^*(\phi(x)\,dx)=\phi(x)\,dx$, where $\phi(x)$ is like (1), (2) or (3) of claim \ref{cl:21}. Then:
\begin{itemize}
\item[i.] If $\phi$ is like in (1) then $h(x)=\d\,x$, where $\d^{m+1}=1$.
\item[ii.] If $\phi$ is like in (2) then $h(x)=\rho\,x$, $\rho\in \C^*$.
\item[iii.] If $\phi$ is like in (3) then $h(x)=\d\,exp(t\,Z)$ for some $t\in \C$, where $\d^{\ell-1}=1$ and $Z=\frac{{x}^\ell}{1+\la\,{x}^{\ell-1}}\,\frac{\pa}{\pa {x}}$. 
\end{itemize}
In particular, in any case, the group $G=\left\{h\in \wh{Diff}(\C,0)\,|\,h^*(\phi(x)\,dx)=\phi(x)\,dx\right\}$ is abelian.
\end{claim}

The proof of claim \ref{cl:22} can be found in \cite{CM} or \cite{LNSSc}.
\vskip.1in
Now, let $p\in D\setminus Sing(\fa^*)$ and $\Si$ be a transverse section to $\fa^*$ with $p\in\Si$. Let $x\in(\C,0)\mapsto \rho(x)\in(\Si,p)$ be a parametrization of $\Si$, so that, we can consider $Diff(\Si,p)=Diff(\C,0)$ and the holonomy group of $D$ in the section as a subgroup $H\sub Diff(\C,0)$.
Let $\wh{x}=x\,u(x)$ be a formal change of variables as in claim \ref{cl:21}: $\eta:=\Pi^*\left(\frac{1}{\wh{f}}\Om\right)|_\Si$ can be written a in (1), (2) or (3) of claim \ref{cl:21}.
\vskip.1in
The definition of holonomy and claim \ref{cl:22} implies that any $h\in H$ satisfies $h^*(\eta)=\eta$. In particular, if $\Om$ has a formal integrating factor then the holonomy group of $D$ must be abelian, contradicting the hypothesis, which proves lemma \ref{l:21}.
\qed
\vskip.1in
Let us continue the proof of theorem \ref{t:2}.
 
\begin{claim}\label{cl:23}
Any $\phi\in \wh{Iso}(\fa)/\wh{Fix}(\fa)$ has a semi-simple representative $\phi_S\in \wh{Iso}(\fa)$. Moreover:
\begin{itemize}
\item[a.] If $\rho\colon\wh{Iso}(\fa)\to GL(2,\C)$ is the group homomorphism $\rho(\phi)=D\phi(0)$ then $ker(\rho)\sub\wh{Fix}(\fa)$.
\item[b.] If $\phi\in\wh{Iso}(\fa)\setminus\wh{Fix}(\fa)$ then $\rho(\phi)$ has finite order in $GL(2,\C)$.
\end{itemize}
\end{claim}

{\it Proof.}
Given $\phi\in\wh{Iso}(\fa)$ let $\phi=\phi_S\circ\phi_U$ be the the Jordan decomposition of $\phi$, where $\phi_S$ is semi-simple and $\phi_U$ unipotent. As we have seen $\phi_S,\phi_U\in\wh{Iso}(\fa)$. By proposition \ref{c:21} there exists a formal nilpotent vector field $Y\in\wh\X_2$ such that $\phi_U=exp(Y)$. On the other hand, by proposition \ref{p:23} $\phi_U\in \wh{Fix}(\fa)$, because otherwise $\Om$ would have a formal integrating factor, which contradicts lemma \ref{l:21}. In particular, the classes of $\phi$ and $\phi_S$ in $\wh{Iso}(\fa)/\wh{Fix}(\fa)$ are the same, which proves the first assertion of the claim.

Assume that $\phi\in \wh{Iso}(\fa)$ and $D\phi(0)=I$ then there exists $Y\in\wh\X_2$ such that $\phi=exp(Y)$. Therefore, as above, $\phi\in\wh{Fix}(\fa)$.
This proves (a).

It remains to prove (b). Let $\phi\in\wh{Iso}(\fa)\setminus\wh{Fix}(\fa)$. By the first assertion we can suppose that $\phi$ has a semi-simple representative in $\wh{Iso}(\fa)$, that we call $\phi_S$. In this case, there exists $\var\in\wh{Diff}(\C^2,0)$ such that $\var^{-1}\circ\phi\circ\var=D\phi(0):=L$ is linear and diagonal in some base of $\C^2$.

Suppose by contradiction that $L$ has infinite order. Let $\ov{H}$ be the Zariski closure in $GL(2,\C)$ of the group $H=\{L^k\,|\,k\in\Z\}$. As we have seen in the proof of theorem \ref{t:1}, $\ov{H}$ is abelian and $dim(\ov{H})\ge1$. Let $\wh{H}=\var\circ \ov{H}\circ\var^{-1}$, so that $\phi_S\in \wh{H}$. Since $\phi_S$ has finite order in $\wh{Iso}(\fa)/\wh{Fix}(\fa)$ there exists $\ell\in\N$ such that $\phi_S^{n\ell}\in \wh{Fix}(\fa)$, $\forall n\in\Z$.
Note that $\ov{H}$ is also the Zariski closure of $\{L^{n\ell}\,|\,n\in\Z\}$. As a consequence, if $\H$ is the Lie algebra of $\ov{H}$ and $\wh\H=\var_*\H$, then for any $Y\in \wh\H$ we have $i_Y\Om=0$ and $exp(Y)\in \wh{Fix}(\fa)$. However, since $L=exp(A)$ for some $A\in\H$, this would imply that $\phi_S=exp(\wh{A})\in \wh{Fix}(\fa)$, where $\wh{A}=\var_*(A)$, a contradiction. 
\qed
\vskip.1in
By the argument of the proof of claim \ref{cl:23} we can construct an injective representation $\rho\colon \wh{Iso}(\fa)/\wh{Fix}(\fa)\to GL(2,\C)$: $\rho(\phi)=D\phi_S(0)$.
Denote $G:=\rho\left(\wh{Iso}(\fa)/\wh{Fix}(\fa)\right)$, the image of $\rho$. Observe that any element $T\in G$ has finite order, say $o(T)\in\N$.
\vskip.1in
The idea is to prove that there exists $k\in\N$ such that $o(T)\le k$, $\forall T\in G$.
It is known that any subgroup of $GL(2,\C)$ with this property is finite and this will finish the proof of theorem \ref{t:2}.
\begin{lemma}\label{l:22}
There exists $k\in\N$ such that $o(T)\le k$, $\forall T\in G$.
\end{lemma}
{\it Proof.}
We will assume first that the resolution process of $\fa$ involves just one blowing-up $\Pi\colon(M,D)\to(\C^2,0)$.
In this case, $D\simeq\p^1$ is not dicritical and contains at least three singularities of $\Pi^*(\fa):=\fa^*$, because the holonomy group of $D$ is not abelian.
Let $Sing(\fa^*)=\{p_1,...,p_k\}$, $k\ge3$, and $\ell=k!$.
We assert that for any $T\in G$ then
\[
T^\ell\in\{\la.\,I\,|\,\la\in\C^*\}\,.
\]
In fact, fix $T\in G$ and $\phi\in \wh{Iso}(\fa)$ semi-simple such that $D\phi(0)=T$.
Assume that $\fa$ is represented by the germ of vector field $X$, with Taylor series $X=\sum_{j\ge n}X_j$, where $X_j$ is homogeneous of degree $j$, $\forall j$, and $X_n\ne0$.
From $\phi^*(X)=u\,X$, $u\in\wh\O_2^*$, we get $T^*(X_n)=u(0)\,X_n$. The singularities $p_1,...p_k$ of $\fa^*$ in $D$ correspond to the $X_n$-invariant directions of $\C^2$. In particular, $T$ induces a permutation of the set $\{p_1,...,p_k\}$ and so $T^\ell(p_j)=p_j$, $1\le j\le k$. A linear isomorphism of $\C^2$ that preserves more than two directions is of the form $\la\,I$, $\la\in\C^*$, which proves the assertion.
\vskip.1in
Let $H=\left<T^\ell\,|\,T\in G\right>\sub\{\la\,I\,|\,\la\in\C^*\}$ and let $S\in H$ with $o(S)=r$.
\begin{claim}\label{cl:24}
\rm We assert that there exists $h\in \wh{Diff}(\C^2,0)$ and $v\in\wh\O_2^*$ such that $h^*(v.\,X):=\wt{X}$ has the Taylor series of the form
\[
\wt{X}=\sum_{m\ge1}\wt{X}_{n+m.r}\,,\,\text{where}\,\,\wt{X}_n=X_n\,.
\]
\end{claim}

In fact, let $\phi\in \wh{Iso}(\fa)$ be semi-simple with $D\phi(0)=S$. 
Let $\phi^*(X)=u.\,X$, where $u\in\wh\O_2^*$. Since $\phi$ is formally linearizable there exists $f\in \wh{Diff}(\C^2,0)$ such that $f^{-1}\circ \phi\circ f=S$.
Let $Y=f^*(X)$, so that $S^*(Y)=\wh{u}.\,Y$, where $\wh{u}=u\circ f$.
Note that $u(0)=\d$ where $\d^r=1$.
This follows from
\[
(S^j)^*(Y)=\wh{u}_j.\,Y\,\,,\,\,\wh{u}_j=\Pi_{i=0}^{j-1}\wh{u}\circ S^i\,.
\]
Since $S^r=I$ we must have $\wh{u}_r=1$ $\implies$ $\d^r=\wh{u}_r(0)=1$.

We assert that there exists $\wh{v}\in \wh\O_2^*$ such that $\wh{v}(0)=1$ and
\begin{equation}\label{eq:6}
S^*(\wh{v}.\,Y)=\wh{u}(0).\wh{v}\,Y\,.
\end{equation}
As the reader can check, the existence of $\wh{v}$ as in (\ref{eq:6}) is equivalent to solve the functional equation
\begin{equation}\label{eq:7}
\wh{v}\circ S=\frac{\wh{u}(0)}{\wh{u}}.\,\wh v:=w.\,\wh{v}\,.
\end{equation}

Note that $\Pi_{j=0}^{r-1}w\circ S^j=1$.
If $\var\in\wh\O_2$ is such that $exp(\var)=w$ then 
\begin{equation}\label{eq:8}
\sum_{j=0}^{r-1}\var\circ S^j=0
\end{equation}
The reader can check that relation (\ref{eq:8}) implies that it is possible to solve the functional equation
\[
\te\circ S-\te=\var\,\,\implies\,\,\wh{v}:=exp(\te)\,\,\text{is a solution of (\ref{eq:7})\,.}
\]
We leave the details for the reader.
\vskip.1in
Set $\wt{X}=\wh{v}.\,Y$. Let $\wt{X}=\sum_{j\ge n}\wt{X}_j$, where $\wt{X}_j$ is homogeneous of degree $j$, $j\ge k$. Note that $\wt{X}_n=X_n$.

We have seen that $S^*(\wt{X})=\d\,\wt{X}$, where $S=\la\,I$, $\la$ is a primitive $r^{th}$ root of the unity and $\d^r=1$.
A direct computation shows that
\[
S^*(\wt{X})=\sum_{j\ge k}\la^{j-1}.\,\wt{X}_j=\sum_{i=0}^{r-1}\la^i.\,Z_i=\d\,\wt{X}=\d\,\sum_{i=0}^{r-1}Z_i\,,
\]
where
\[
Z_i=\sum_{m=0}^\infty\wt{X}_{i+1+m.r}\,,\,1\le i\le r-1\,.
\]
Therefore,
\[
\sum_{i=0}^{r-1}(\la^i-\d)Z_i\equiv0\,\,.
\]
Since $Z_{n-1}=X_n+h.o.t.\ne0$ we have $\d=\la^{n-1}$ and $Z_i=0$ if $i\ne n-1$.
In particular, we get
\[
\wt{X}=Z_{n-1}=X_n+\sum_{m\ge1}\wt{X}_{n+m.r}\,,
\]
which proves claim \ref{cl:24}.
\qed
\vskip.1in
Let $\wt{X}_j=P_j\frac{\pa}{\pa x}+Q_j\frac{\pa}{\pa y}$, where $P_j$ and $Q_j$ are homogeneous of degree $j$. 
In the chart $\Pi(t,x)=(x,t\,x)=(x,y)$ of the blow-up $\Pi$ we obtain
\[
\Pi^*(X_j)=x^{j-1}\left(f_j(t)\frac{\pa}{\pa t}+x\,P_j(1,t)\frac{\pa}{\pa x}\right)\,,\,\text{where}\,\,f_j(t)=Q_j(1,t)-t\,P_j(1,t)\,.
\]
This implies $\Pi^*(\wt{X})=x^{n-1}\left[x\,G(t,x^r)\frac{\pa}{\pa x}+F(t,x^r)\frac{\pa}{\pa t}\right]$, where
\[
\left\{
\begin{matrix}
G(t,x^r)=\sum_{m\ge0}f_{n+mr}(t)\,x^{m\,r}\\
F(t,x^r)=\sum_{m\ge0}P_{n+mr}(1,t)\,x^{mr}\\
\end{matrix}
\right..
\]
The strict transform of $\Pi^*(\fa)$ is therefore defined by the form $\Om^*:=x\,G(t,x^r)\,dt-F(t,x^r)\,dx$.

Let $\g\colon [0,1]\to D\setminus\{p_1,...,p_k\}$ be a $C^1$ closed curve with $\g(0)=\g(1)=t_o\in D$. The holonomy $h_\g$ of the curve $\g$ calculated in the section $(t=t_o)$ is the solution of the differential equation
\begin{equation}\label{eq:9}
\frac{dx}{ds}=x.\,\frac{G(\g(s),x^r)}{F(\g(s),x^r)}\,\g'(s)
\end{equation}
with initial condition $x(0)=x$. If $X(s,x)$ is this solution then $h_\g(x)=X(1,x)$.

Note that $h_\g$ is necessarily of the form
\[
h_\g(x)=\la\,x\,(1+H_\g(x^r))\,.
\]
In fact, if we consider the ramification $y=x^r$ applyed in (\ref{eq:9}) we obtain
\[
\frac{dy}{ds}=r\,x^{r-1}\frac{dx}{ds}=r\,x^r\,\frac{G(\g(s),x^r)}{F(\g(s),x^r)}\,\g'(s)=r\,y\,\frac{G(\g(s),y)}{F(\g(s),y)}\,\g'(s)\,.
\]
If $Y(s,y)$ is the solution of the above equation with initial condition $Y(0,y)=y$ then
\begin{itemize}
\item[(1).] $Y(s,y)=y.\,U(s,y)$ where
\[
U(s,0)=exp\left(r\int_0^s\frac{G(\g(s),0)}{F(\g(s),0)}\,\g'(s)\,ds\right)\,.
\]
\item[(2).] $X(s,x)=x\,(U(s,x^m))^{1/r}$. This implies the assertion.
\end{itemize}
Note also that
\[
h_\g'(0)=exp\left(\int_\g\frac{G(z,0)}{F(z,0)}\,dz\right)\,.
\]

Since the holonomy of $D$ is non abelian, there are $\a,\be\in \Pi_1(D\setminus\{p_1,...,p_k\},t_o)$ such that if $\g=[\a,\be]:=\a.\be.\a^{-1}.\be^{-1}$, then $h_\g\ne id$, but is tangent to the identity.
Therefore, the formal diffeomorphism is necessarily of the form
\[
h_\g(x)=x+a.\,x^{m.\,r}+h.o.t.\,\,,\,\,\text{where $a\ne0$ and $m\ge1$}.
\]
The integer $m.\,r$ is the order of tangency of $h_\g$ with the identity. It is a formal invariant of $h_\g$, in the sense that if $g=f^{-1}\circ h_\g\circ f$, where $f\in \wh{Diff}(\C,0)$ then the order of tangency of $g$ with the identity is also $m.\,r$.

Suppose by contradiction that the assertion of lemma \ref{l:22} is false. This implies that the set $\{o(T)\,|\,T\in G\}$ is unbounded, so that there is $T\in G$ such that $o(T)=r'>m.\,r$. But, as we have seen above, this implies that the order of tangency of $h_\g$ with the identity is a multiple $m'.\,r'>m.\,r$, a contradiction.
Therefore, the set $\{o(T)\,|\,T\in G\}$ is finite, as we wished.
\vskip.1in
Suppose now that the resolution process $\Pi\colon(M,E)\to(\C^2,0)$ involves more than one blowing-up.
Let $\wh{Diff}(M,E)$ be the set of germs at $E$ of formal diffeomorphims of $M$. 

We assert that for any $\phi\in\wh{Iso}(\fa)$ there is $\wt\phi\in \wh{Diff}(M,E)$ such that the diagram below commutes
\begin{equation}\label{eq:lift}
\begin{matrix}
(M,E)&\overset{\wt\phi}\longrightarrow&(M,E)\\
\Pi\downarrow&\,&\downarrow\Pi\\
(\C^2,0)&\underset{\phi}\longrightarrow&(\C^2,0)\\
\end{matrix}
\end{equation}

The proof is done following the resolution process. To give an idea, we will do the two first steps.
Since $\phi(0)=0$, when we perform the first blow-up $\Pi_1\colon(M_1,D_1)\to(\C^2,0)$ then there is $\phi_1\in \wh{Diff}(M_1,D_1)$ that lifts $\phi$: $\Pi_1\circ\phi_1=\phi\circ\Pi_1$.

The formal diffeomorphism $\phi_1$ converges in the divisor and $\phi_1|_{D_1}$ is an automorphism of $D_1$. Moreover, it preserves the germ of the strict transform of $\fa$, say $\fa_1$, along $D_1$. In particular, it induces a permutation in the set of singularities $Sing(\fa_1)\sub D_1$.

Let $p_1\in Sing(\fa_1)$ be a non-simple singularity of $\fa_1$. Its orbit by $\phi_1$ is periodic: $p_2=\phi_1(p_1)$, $p_3=\phi_1(p_2)$, ..., $p_1=\phi_1(p_{\ell})$.
Since $\phi_1^*(\fa_1)=\fa_1$ the germs of $\fa_1$ at all $p_{j's}$ are equivalent and to continue the blowing-up process we have to blow-up in all these singularities, thus obtaining another step of the process $\Pi_2\colon (M_2,E_2)\to(M_1,D_1)$, with $E_2=\wt{D}_1\cup D_1'\cup...\cup D_\ell'$, where $\wt{D}_1$ is the strict transform of $D_1$ and $D_j'=\Pi_2^{-1}(p_j)$.

In this case, we can lift $\phi_1$ to a formal germ of diffeomorphism $\phi_2\in\wh{Diff}(M_2,E_2)$ such that
\begin{itemize}
\item[(3).] $\Pi_2\circ\phi_2=\phi_1\circ\Pi_2$.
\item[(4).] $\phi_2|_{E_2}$ is an automorphism of $E_2$.
\item[(5).] $\phi_2(\wt{D}_1)=\wt{D}_1$ and $\phi_2(D_j')=D_{j+1}'$, $\forall 1\le j\le \ell-1$.
\item[(6).] If $\fa_2=\Pi_2^*(\fa_1)$ then $\phi_2^*(\fa_2)=\fa_2$.
\end{itemize}

Continuing this process inductively, at the end we will find a germ of formal diffeomorphism $\wt\phi\in\wh{Diff}(M,E)$ that makes the diagram (\ref{eq:lift}) to commute.
Moreover, if $\wt\fa=\Pi^*(\fa)$ then $\wt\phi^*(\wt\fa)=\wt\fa$.
The restriction $\wt\phi|_E$ is an automorphism of $E$ and permutes the irreducible components of $E$ and also $\wt\phi(Sing(\wt\fa))=Sing(\wt\fa)$, $\wt\fa=\Pi^*(\fa)$. 
\vskip.1in
Let $D$ be an irreducible component of $E$ with non abelian holonomy. Since $E$ has a finite number $k$ of irreducible components, if $m=k!$ then for any $\phi\in\wh{Iso}(\fa)$ we have ${\wt\phi}^m(D)=D$, where $\wt\phi$ lifts $\phi$. Let $\#\,(Sing(\wt\fa)\cap D)=k'$. Since the holonomy of $D$ is non abelian we have $k'\ge3$. In particular, if $m'=k!.k'!$ then $\wt\phi^{m'}|_D=id_D$, the identity of $D$.

Assume that the divisor $D$ was obtained blowing-up a singularity $q$ in a previous step of the process. Taking local coordinates at $q$ we can assume that the blowing-up is
$\Pi\colon(\wt M,D)\to(\C^2,q=0)$ and the germ at $q$ of the previous foliation in the process was
$\fa'$, so that $\Pi^*(\fa')$ is the germ of $\wt\fa$ along $D$.
The group $\wt{G}:=\left<\wt\phi^{m'}|_{M'}\,|\,\phi\in \wh{Iso}(\fa)\right>\sub\wh{Diff}(M',D)$ is a subgroup of $\wh{Iso}(\wt\fa|_{M'})$.

Repeating the previous argument it can be shown that $\wh{Iso}(\wt\fa|_{M'})/\wh{Fix}(\wt\fa|_{M'})$ is finite. In particular, there exists $k\in N$ such that for any semi-simple $\phi\in \wh{Iso}(\fa)\setminus\wh{Fix}(\fa)$ we have $\wt\phi^k|_{M'}=id_{M'}$, which implies that $\wt\phi^k=id_M$ and $\phi^k=id_{\C^2}$, proving lemma \ref{l:22} and theorem \ref{t:2}.
\qed

\subsection{Examples with formal, but without meromorphic integrating factor}\label{ss:24}
The purpose of this section is to prove the following result:

\begin{thm}\label{t:3}
There exist germs of foliations $\fa$ on $(\C^2,0)$ with the following properties:
\begin{itemize}
\item[(a).] $Iso(\fa)/Fix(\fa)$ is non-countable and contains elements of infinite order.
\item[(b).] $\fa$ is not liouville integrable: it cannot be defined by a closed meromorphic 1-form.
\end{itemize}
\end{thm}

{\it Proof.}
The proof is based in \cite{pm} and in the fact that any germ $h\in Diff(\C,0)$ can be realized as the holonomy of a separatrix of a germ holomorphic vector field $X\in\X(\C^2,0)$ such that $DX(0)$ is linear and diagonal:

\begin{claim}\label{cl:25}
\rm Let $h\in Diff(\C,0)$ be a germ of biholomorphism with $h'(0)=\la\in\C^*$. Then there exists a holomorphic vector field $X\in\X(\De)$, $\De=\{(x,y)\in\C^2\,;\,|x|<2\,,\,|y|<\ep\}$ of the form
\[
X(x,y)=x\frac{\pa}{\pa x}+y\,b_h(x,y)\frac{\pa}{\pa y}\,,
\]
where $\la=e^{2\pi i\a_h}$, $\a_h=b_h(0,0)$, and the holonomy of the curve $\be(t)=(e^{2\pi it},0)$, $t\in[0,1]$, contained in the leaf $(y=0)\setminus\{(0,0)\}$, in the transversal section $\Si:=(x=1)$ is $y\mapsto h(y)$.
\end{claim}
{\it Proof.} 
When $h$ is linearizable theorem \ref{t:3} is immediate. When $h$ is non-linearizable then $|\la|=1$. When $\la$ is a root of unity the proof can be found in \cite{mr}, whereas when $\la$ is not a root of unity the proof can be found in \cite{pmy}.
\qed
\vskip.1in
Let $\fa=\fa_X$ be the germ of foliation defined by the vector field $X$.
The foliation $\fa$ satisfies the following properties:
\begin{itemize}
\item[(1).] Outside the axis $(x=0)$ it is transverse to the vertical fibration $(x=ct)$.
\item[(2).] The leaf $L_y$ of $\fa$ through the point $(1,y)$ cuts the fiber $(x=1)$ exactly at the points of the form $(1,\wt{y})$, where $\wt{y}$ belongs to the pseudo-orbit of $h$: $\wt{y}=h^n(y)$, $n\in \Z$.
\end{itemize}
More precisely, given $n\in\Z$ let $Dom(h^n)=$ the connected component of $0\in\C$ of the set $\{y\in\C\,|\,|y|<\ep$ and $h^j(y)$ is defined for all $j$ with $0\le|n-j|\le |n|\}$. 
If $y\in Dom(h^n)$ then $(1,h^j(y))\in L_y\cap(x=1)$ for all $j$ with $0\le|n-j|\le |n|$.

We will assume also that $\a_h\notin \R_+$. With this condition, the saturation by $\fa$ of the set $\Si_\ep=\{(1,y)\,|\,|y|<\ep\}$
\[
Sat(\Si_\ep)=\bigcup_{y\in \Si_\ep}L_y
\]
contains a set of the form $D_\d\times D_\ep\setminus \{x=0\}$, where $D_r$ denotes the disc $\{z\in\C\,|\,|z|<r\}$.
\vskip.1in
Given $h_1\in Diff(\C,0)$ a germ of diffeomorphism commuting with $h$, we will construct a germ $\Phi_{h_1}\in Iso(\fa)\sub Diff(\C^2,0)$.
The construction will be done in such a way that:
\begin{itemize}
\item[(3).] $\Phi_{h_1}(x,y)=(x,f(x,y))$, so that $\Phi_{h_1}$ preserves the fibers $(x=ct)$.
\item[(4).] $\Phi_{h_1}(1,y)=(1,h_1(y))$\, ($f(1,y)=h_1(y)$).
\end{itemize}

\begin{rem}\label{r:22}
{\rm We will see that there is only one $\Phi_{h_1}\in Iso(\fa)$ satisfying (3) and (4).}
\end{rem}

In order to formalize the construction of $\Phi_{h_1}$ we consider the universal covering of
$\De\setminus\{x=0\}$, $\Pi\colon B\times D_\ep\to \De\setminus\{x=0\}$,
\[
\Pi(z,y)=(e^z,y)\,,
\]
where $B=\{z\in\C\,|\,Re(z)\in(-\infty,\ell og(2))\}$ and $\De=D_2\times D_\ep$.

The pull-back foliation $\wt\fa:=\Pi^*(\fa)$ is defined in $B\times D_\ep$ by the vector field
\[
\wt{X}(z,y)=\Pi^*(X)(z,y)=\frac{\pa}{\pa z}+y.\,b_h(e^z,y)\frac{\pa}{\pa y}\,.
\]
Note that this foliation has no holonomy, in the sense that if $\wt L$ is a leaf of $\wt\fa$ then it cuts any transversal $(z=ct)$ in at most one point.
Denote by $\wt L_y$ the leaf of $\wt\fa$ such that $\wt L_y\cap (z=0)=(0,y)$.

Note that, by the definition of holonomy (of $\fa$), if $y\in Dom(h^n)$ then $\wt L_y\cap(z=2\,n\,\pi\,i)=(2\,n\,\pi\,i,h^n(y))$. 

Now, we define the covering $\wt\Phi_{h_1}$ of (the future) $\Phi_{h_1}$ as the germ along $(y=0)=B\times\{0\}$ of map that satisfies:
\begin{itemize}
\item[(5).] $\wt\Phi_{h_1}$ preserves the fibers $(z=ct)$: $\wt\Phi_{h_1}(z=z_o)\sub(z=z_o)$. 
\item[(6).] $\wt\Phi_{h_1}(\wt L_y)=\wt L_{h_1(y)}$.
\end{itemize}
Since $\wt\fa$ has no holonomy (5) and (6) define an unique germ of holomorphic diffeomorphism $\wt\Phi_{h_1}$ along $B\times\{0\}$ such that $\wt\Phi_{h_1}(B\times\{0\})=B\times\{0\}$.

Now, we will see that there exists a germ $\Phi_{h_1}\in Iso(\fa)\sub Diff(\C^2,0)$ such that $\Phi_{h_1}\circ\Pi=\Pi\circ\wt\Phi_{h_1}$.
By (5) and (6), the extension of $\wt\Phi_{h_1}$ to the fiber $(t=t_o)$ is done using the holonomy of $\wt\fa$, say $H_{t_o}\colon(t=t_o)\to(t=0)$, so that
\begin{equation}\label{eq:12}
\wt\Phi_{h_1}(t_o,y)=(t_o,H_{t_o}^{-1}(h_1(H_{t_o}(y))))\,.
\end{equation}
For instance, $H_{2\pi i}(y)=h^{-1}(y)$ and so
\[
\wt\Phi_{h_1}(2\pi i,y)=(2\pi i,h(h_1(h^{-1}(y))))=(2\pi i,h_1(y))\,,
\]
because $h$ and $h_1$ commute. In particular,
\[
\Pi\circ\wt\Phi_{h_1}(2\pi i,y)=(1,h_1(y))=\Pi(0,h_1(y))=\Pi\circ\wt\Phi_{h_1}(0,y)\,.
\]
Similarly, $\Pi\circ\wt\Phi_{h_1}(2k\pi i,y)=\Pi\circ\wt\Phi_{h_1}(0,y)$ for all $k\in\Z$ and $y\in Dom(h^k)$.
This implies that we can define $\Phi_{h_1}(1,y)=(1,h_1(y))$ for all $y\in Dom(h_1)$.

The extension of $\Phi_{h_1}$ to the fibers $(x=x_o)$, $x_o\ne0$, can be done by using (\ref{eq:12}). We leave the details to the reader.
It remains to prove that $\Phi_{h_1}$ can be extended to the fiber $(x=0)$.

When $\a=\a_h\notin \R$ then, by Poincar\'es linearization theorem, we can assume that $X$ is linear. In this case, $h(y)=\la\,y$, where $\la=exp(2\pi i\a)$. If $h_1$ commutes with $h$ then $h_1$ is also linear and $\Phi_{h_1}(x,y)=(x,h_1(y))$, as the reader can check. 

When $\a\in\R_-$ and $h$ is non-linearizable, then we have to use that the saturation of the transversal $(x=1)$ by $\fa$ contains a set of the form
$D_\d\times D_\ep\setminus(x=0)$. In this case, it can be proved that $\Phi_{h_1}$ is bounded in the set $D_\d\times D_\ep\setminus(x=0)$ and so $\Phi_{h_1}$ can be extended to $\{0\}\times D_\ep$ by Riemann's extension theorem (see \cite{gr}). We leave the details to the reader.

Now, we use a construction of Perez Marco. In \cite{pm} he proves the following result:
\vskip.1in
{\bf Theorem.}
{\it There exists non-linearizable germs of diffeomorphisms $h\in Diff(\C,0)$ of the form $h(y)=\la\,y+h.o.t.$, with $|\la|=1$ and not a root of unity, whoose centalizer
\[
C(h)=\{g\in Diff(\C,0)\,|\,g\circ h\circ g^{-1}\circ h^{-1}=Id\}
\]
is a Cantor set, in the sense that the set $\{g'(0)\,|\,g\in C(h)\}$ is a Cantor set of $S^1=\{\mu\in\C\,|\,|\mu|=1\}$. In particular, $C(h)$ is non countable. In fact, in $C(h)$ there are infinitely many elements of finite order and a non countable set of elements of infinite order.}

Now, we take a vector field $X$ like in claim \ref{cl:25}, associated to $h$ like in Perez Marco's theorem. The foliation $\fa_X$, associated to $X$, cannot be defined by a closed meromorphic 1-form, and $Iso(\fa_X)/Fix(\fa_X)$ is non-countable.

This finishes the proof of theorem \ref{t:3}.
\qed

\section{Conical and logaritmic foliations}\label{ss:3}
In this section we study $Iso(\fa)$ when $\fa$ is a conical or a logarithmic foliation.
\vskip.1in
\subsection{Preliminaries and statement of the results.}\label{ss:31}
A conical foliation $\fa$ on $\C^n$ induces a foliation of the same codimension, say $\wt\fa$, on the  projective space $\p^{n-1}$. Denote by $Aut(\wt\fa)$ the subgroup of $Aut(\p^{n-1}_K)$, of automorphisms of $\p^{n-1}_K$ preserving $\wt\fa$:
\[
Aut(\wt\fa)=\left\{\Phi\in Aut(\p^{n-1}_K)\,|\,\Phi^*(\wt\fa)=\wt\fa\right\}\,.
\]
We have the following:

\begin{thm}\label{t:4}
If $\fa$ is conical and $\wt\fa$ is the foliation induced by $\fa$ on the projective space $\p^{n-1}$ then $Iso(\fa)/Fix(\fa)$ and $\wh{Iso}(\fa)/\wh{Fix}(\fa)$ are isomorphic to $Aut(\wt\fa)$.
\end{thm}

When $\fa$ is a conical homogeneous foliation of dimension one and the degree of the foliation $\wt\fa$, as a foliation of $\p^{n-1}$, is $d$ then the degree of the homogeneous form $\Om$ on $\C^n$ inducing $\fa$ is $d+1$ (see \cite{LW}). Denote by $\fol(d,n-1)$ the set of 1-dimensional foliations on
$\p^{n-1}$ of degree $d$. It is known that $\fol(d,n-1)$ can be identified with a Zariski open and dense subset of some $\p^N$.
As a consequence of theorem \ref{t:4} we have the following:

\begin{Cor}\label{c:1}
If $d\ge2$ then $\fol(d,n-1)$ contains a Zariski open and dense subset, say $G$, such that for any $\G\in G$ we have $Iso(\Pi^*(\G))=Fix(\Pi^*(\G))$, where $\Pi\colon\C^n\setminus\{0\}\to\p^{n-1}$ is the canonical projection.
\end{Cor}

When $\Om\in\Om^1(\C^n)$ is homogeneous and non-conical, then it has a holomorphic integrating factor: if $f:=i_R\Om\ne0$ then $f^{-1}.\Om$ is closed. In particular, if the decomposition of $f$ into irreducible factors is $f_1^{k_1}...f_r^{k_r}$ then
\[
\frac{\Om}{f}=\sum_{j=1}^r\la_j\frac{df_j}{f_j}+d\left(\frac{H}{f_1^{k_1-1}...f_r^{k_r-1}}\right)\,,
\] 
where $\la_j\in\C$, $1\le j\le r$, and $f_1$,...,$f_r$,$H$ are all quasi-homogeneous with respect to $S$. When $f=f_1...f_r$, is reduced, then in the above formula we have $H=0$ and $\la_1...\la_r\ne0$. In this case, the form $f^{-1}.\,\Om$ is logarithmic. 

In section \ref{ss:33} we will consider germs of closed logarithmic 1-forms in general, that is when the divisor of poles is reduced and not necessarily quasi-homogeneous. In this case, if $\Om$ has pole divisor $F=f_1...f_r$ then
\[
\Om=\sum_{j=1}^r\la_j\frac{df_j}{f_j}+dh\,,
\]
where $\la_1,...,\la_r\in\C^*$ and $h\in\O_n$.

\begin{rem}
\rm Multiplying $f_1$ by an unity $u\in\O_n^*$ we can suppose that $dh=0$.
In fact, since $\la_1\ne0$, if we set $\wt{f}_1=u\,f_1$, where $u=exp(-h/\la_1)$, then
\[
\Om=\la_1\frac{d\wt f_1}{\wt f_1}+\sum_{j\ge2}\la_j\frac{df_j}{f_j}\,.
\]

Other remarks about a 1-form $\Om$ as above are the following

\begin{itemize}
\item[a.] $\Om$ has a holomorphic (or formal) non-constant first integral if, and only if, there exists $\la\in\C^*$ such that $(\la_1/\la,...,\la_n/\la)\in\N^n$.
\item[b.] $\Om$ has a meromorphic (or formal meromorphic) non-constant first integral if, and only if, there exists $\la\in\C^*$ such that $(\la_1/\la,...,\la_n/\la)\in\Z^n$.
\end{itemize}
\end{rem}

Another fact is the following:

\begin{prop}\label{p:4}
Let $\Om$ be a germ at $0\in\C^n$ of closed meromorphic 1-form, without non-constant meromorphic first integral. If $\Phi\in \wh{Iso}(\fa_\Om)$ then $\Phi^*(\Om)=\d.\,\Om$, where $\d$ is a root of unity. 
\end{prop}

\begin{ex}
\rm In general $\d=1$ in proposition \ref{p:4}, as we will see in the proof. Let us see an example in which $\d=e^{2\pi i/3}$: we set $\Om=\frac{dx}{x}+\d\frac{dy}{y}+\d^2\frac{dz}{z}$ and $\Phi(x,y,z)=(z,x,y)$. As the reader can check we have $\Phi^*(\Om)=\d.\,\Om$.
\end{ex} 

In the next result we give conditions implying that a logarithmic 1-form is holomorphically equivalent to a homogeneous form.
\begin{thm}\label{t:5}
Let $\Om=\sum_{j=1}^r\la_j\frac{df_j}{f_j}$, where $f_1,...,f_r$ are irreducible and $\fa_\Om$ has no meromorphic first integral. Suppose also that:
\begin{itemize}
\item[(a).] There exists $\Phi\in \wh{Iso}(\fa_\Om)$ such that $D\Phi(0)=\rho.\,I$ and,
either $\rho$ is not a root of unity, or the class of $\Phi$ in $\wh{Iso}(\fa_\Om)/\wh{Fix}(\fa_\Om)$ has infinite order.
\item[(b).] The first non zero jet of $f_j$ is $J^{k_j}_0(f_j)=h_j$ and $\sum_{j=1}^rk_j.\,\la_j\ne0$.
\end{itemize}
Then there exists $\phi\in Diff(\C^n,0)$ such that
\begin{equation}\label{eq:} 
\phi^*(\Om)=\sum_{j=1}^r\la_j\frac{dh_j}{h_j}\,.
\end{equation}
In particular, $\fa_\Om$ is holomorphically equivalent to a homogeneous foliation.
\end{thm}

As a consequence we have:

\begin{Cor}\label{c:2}
Let $\Om=\sum_{j=1}^r\la_j\frac{df_j}{f_j}$, where $f_1,...,f_r\in\O_n$ are irreducible.
Suppose further that:
\begin{itemize}
\item[(a).] $\fa_\Om$ has no meromorphic first integral.
\item[(b).] The first non-zero jet of $f_j$ is $J^{k_j}_0f_j=h_j$ and $\sum_{j=1}^rk_j\,\la_j\ne0$.
\item[(c).] There exists a formal diffeomorphism $\wh\psi\in \wh{Diff}(\C^n,0)$ such that
$\wh\psi^*(\Om)=\Om_h$, where $\Om_h=\sum_{j=1}^r\la_j\frac{dh_j}{h_j}$.
\end{itemize}
Then there exists $\psi\in Diff(\C^n,0)$ such that $\psi^*(\Om)=\Om_h$.
\end{Cor}

The proof of corollary \ref{c:2} is based in the fact that if $T=\rho.\,I$ then $T\in Iso(\fa_{\Om_h})$.
This of course, implies that hypothesis (a) of theorem \ref{t:5}.
\vskip.1in
We state below a condition implying that for any $\Phi\in Iso(\fa_\Om)$ there exists $N\in\N$ such that $D\Phi^N(0)=\rho.\,I$, where $\rho\in\C^*$.
\begin{defi}\label{d:3}
\rm Let $H=(h_1,...,h_\ell)$ be a set of homogeneous polynomials, not necessarily of the same degree. Set
\[
\I(H)=\{T\in GL(n,\C)\,|\,h_j\circ T=\a_j.\,h_j\,\,,\,\a_j\in\C^*\,,\,1\le j\le \ell\}\,.
\]
Note that:
\begin{itemize}
\item[I.] $\I(H)$ is a closed sub-group of $GL(n,\C)$.
\item[II.] $\I(H)\sup\C.I=\{\rho.I\,|\,\rho\in\C^*\}$, where $I$ is the identity in $GL(n,\C)$.
\end{itemize}
We say that $H=(h_1,...,h_\ell)$ is {\it rigid} if $\I(H)=\C.\,I$.
\end{defi}

\begin{ex}
\rm If $H=(z_1,...,z_n)\sub\O_n$ then $\I(H)=$
\[
\{T\,|\,T(z_1,...,z_n)=(\la_1.\,z_1,...,\la_n.\,z_n)\,\,,\,\,\la_1,...,\la_n\in\C^*\}\,.
\]
Two examples of rigid sets are $H_1=(z_1,...,z_n,z_1+z_2+...+z_n)$ and $H_2=(z_1.\,z_2^2...z_n^n,z_1+...+z_n)$.
\end{ex}

\begin{Cor}\label{c:3}
Let $\Om=\sum_{j=1}^r\la_j\frac{df_j}{f_j}$ and assume that:
\begin{itemize}
\item[(a).] $f_1,...,f_r$ are irreducible and $\fa_\Om$ has no meromorphic first integral.
\item[(b).] The first non zero jet of $f_j$ at $0\in\C^n$ is
$J_0^{k_j}(f_j):=h_j$, where $k_j\ge1$, and the set $H:=(h_1,...,h_r)$ is rigid.
\item[(c).] $\sum_{j=1}^rk_j.\,\la_j\ne0$.
\item[(d).] There exists $\Phi\in \wh{Iso}(\fa_\Om)$ whose class in
$\wh{Iso}(\fa_\Om)/\wh{Fix}(\fa_\Om)$ has infinite order.
\end{itemize}
Then there exists $\phi\in Diff(\C^n,0)$ such that $\phi^*(\Om)=\sum_{j=1}^r\la_j\frac{dh_j}{h_j}$.
\end{Cor}

Another consequence of theorem \ref{t:5} is the following:
\begin{Cor}\label{c:4}
Let $\Om=\sum_{j=1}^r\la_j\frac{df_j}{f_j}$ and assume that:
\begin{itemize}
\item[(a).] $f_1,...,f_r$ are irreducible and $\fa_\Om$ has no meromorphic first integral.
\item[(b).] The first non zero jet of $f_j$ at $0\in\C^n$ is
$J_0^{k_j}(f_j):=h_j$, where $k_j\ge1$, and the set $H:=(h_1,...,h_r)$ is rigid.
\item[(c).] $h_1$ is not irreducible.
\item[(d).] $\sum_{j=1}^rk_j.\,\la_j\ne0$.
\end{itemize}
Then $Iso(\fa_\Om)/Fix(\fa_\Om)\simeq\wh{Iso}(\fa_\Om)/\wh{Fix}(\fa_\Om)$ and both are isomorphic to the same finite sub-group of $GL(n,\C)$.
\end{Cor}

\begin{ex}
\rm An example for wich $Iso(\fa_\Om)=Fix(\fa_\Om)$ and $\wh{Iso}(\fa_\Om)=\wh{Fix}(\fa_\Om)$
is
$\Om=\la_1\frac{df_1}{f_1}+\la_2\frac{df_2}{f_2}+\la_3\frac{df_3}{f_3}+\la_4\frac{df_4}{f_4}$, where $\la_j/\la_i\notin \Q$, $\forall$ $i\ne j$, $f_1=z_1^2+z_2^3$, $f_2=z_1$, $f_3=z_3$ and $f_4=z_1+z_2$.

Another example was done in example \ref{ex:c4}.
\end{ex}

In the case of dimension two we have the following:
\begin{Cor}\label{c:5}
Let $\Om=\sum_{j=1}^r\la_j\frac{df_j}{f_j}\in\Om^1(\C^2,0)$, where $f_1,...,f_r\in\O_2$ are all irreducible and relatively primes two by two.
Suppose that:
\begin{itemize}
\item[(a).] $\fa_\Om$ has no meromorphic first integral.
\item[(b).] The first non-zero jet of $f_j$ is $J^{k_j}_0(f_j)=h_j$, $1\le j\le r$, and the set $(h_1,...,h_r)$ is rigid.
\item[(c).] $\sum_{j=1}^rk_j\la_j\ne0$.
\item[(d).] There exists $\Phi\in\wh{Iso}(\fa_\Om)$ whoose class in
$\wh{Iso}(\fa_\Om)/\wh{Fix}(\fa_\Om)$ has infinite order.
\end{itemize}
Then $k_1=...=k_r=1$, $h_1,...,h_r$ are linear forms and there coordinates $(z_1,z_2)$ such that $\Om$ is holomorphically equivalent to
\begin{equation}\label{eq:d2}
\la_1\frac{dz_1}{z_1}+\la_2\frac{dz_2}{z_2}+\sum_{j=3}^r\la_j\frac{dh_j}{h_j}\,.
\end{equation}
\end{Cor}

Corollary \ref{c:5} is a direct consequence of corollary \ref{c:3} and the fact that any homogeneous polynomial in $\C[z_1,z_2]$ can be decomposed into linear factors. We leave the details to the reader.
 
\subsection{Conical foliations: proof of theorem \ref{t:4}}\label{ss:32}

Let $\om$ be a conic integrable p-form on $\C^n$.
We will assume that $cod(Sing(\om))\ge2$. The form $\om$ is homogeneous and defines a codimension p foliation $\fa_\om$ on $\C^n$ and a foliation $\wt\fa_\om$ on $\p^{n-1}$.
We want to prove that $Aut(\wt\fa_\om)\simeq Iso(\fa_\om)/Fix(\fa_\om)$.

Given $\Phi\in Iso(\fa_\om)$ we denote as $\ov\Phi$ its class in $Iso(\fa_\om)/Fix(\fa_\om)$.
\begin{lemma}\label{l:31}
Let $\Phi\in \wh{Iso}(\fa_\om)$ and $\Phi_1=D\Phi(0)$ be the linear part of $\Phi$ at $0$.
Then:
\begin{itemize}
\item[(a).] $\Phi_1\in Iso(\fa_\om)$.
\item[(b).] $\Phi$ and $\Phi_1$ define the same class at
$\wh{Iso}(\fa_\om)/\wh{Fix}(\fa_\om)$.
\end{itemize}
\end{lemma}

{\it Proof.}
The condition $\Phi\in \wh{Iso}(\fa_\om)$ is equivalent to $\Phi^*(\om)=u.\,\om$, where $u\in\wh\O_n^*$. Since the coefficients of $\om$ are homogeneous of the same degree, the first non-zero jet of $u.\,\om$ is
$u(0).\,\om$, whereas the first non-zero jet of $\Phi^*(\om)$ is $\Phi_1^*(\om)$.
Hence $\Phi_1^*(\om)=u(0).\,\om$ and $\Phi_1\in Iso(\fa_\om)$.

In particular, we can write $\Phi=\Phi_1\circ \var$, where $\var\in \wh{Iso}(\fa_\om)$ and $D\var(0)=I$.
We assert that $\var\in \wh{Fix}(\fa_\om)$.

To see this, let us blow-up once at the origin $0\in\C^n$. If we denote this blow-up by
$\Pi\colon(\wt\C^n,0)\to(\C^n,0)$ then the exceptional divisor $E=\Pi^{-1}(0)\simeq\p^{n-1}$ and there exists a germ
$\wh\var\in\wh{Diff}(\C^n,0)$ such that $\Pi\circ \wh\var=\var\circ \Pi$.

Recall also that $\wt\C^n$ is a linear bundle over $E$, say $P\colon\wt\C^n\to E$, where $E$ is the zero section and the fibers
$P^{-1}(a)$, $a\in E$, project by $\Pi$ onto the lines of $\C^n$ through the origin.

Let $\wh\fa_\om=\Pi^*(\fa_\om)$. Since $\fa_\om$ is conical with respect to the radial vector field
$R=\sum_jz_j\frac{\pa}{\pa z_j}$, the lines through the origin of $\C^n$ are contained in the leaves or in the singular set of $\fa_\om$, so that the fibers of $P$ are contained in the leaves or in the singular set of $\wh\fa_\om$. In particular the leaves and the singular set of $\wh\fa_\om$ are transverse to $E$. 

Finally, we note that $\wt\fa_\om=i_E^*(\wh\fa_\om)$, where $i_E\colon E\to\wh\C^n$ is the inclusion.  In other words, $\wt\fa_\om$ can be viewed as the intersection of $\wh\fa_\om$ with the zero section $E$.
Since $d\var(0)=I$ we get $\wh\var|_E=id_E$, the identity map of $E$, which implies that $\wh\var$ preserves the leaves and the singular set of $\wh\fa_\om$. This implies that $\var\in Fix(\fa_\om)$.
\qed
 
In particular, any $\ov\Phi\in \wh{Iso}(\fa)/\wh{Fix}(\fa)$ has a linear representative $\Phi_1\in GL(n,\C)$.

On the other hand, any $T\in GL(n,\C)$ induces an automorphism $\wt T\in Aut(\p^{n-1})$.
Let $Q\colon GL(n,\C)\to Aut(\p^{n-1})$ be the group homomorphism given by $Q(T)=\wt T$, where $\wt T$ is as before.
Recall that $Q$  is surjective and 
\begin{equation}\label{eq:14}
ker(Q)=\C^*.\,I=\{\rho.\,I\,|\,\rho\in\C^*\}\,.
\end{equation}

Finally, we have seen in lemma \ref{l:31} that $Q(Iso(\fa))\sub Aut(\wt\fa)$ and that $\Phi_1\in Fix(\fa)$ if, and only if, $Q(\Phi_1)=id$, where $id=$ identity of $Aut(\p^{n-1})$.
Therefore (\ref{eq:14}) implies that
\[
Aut(\wt\fa)\simeq \wh{Iso}(\fa)/ker(Q|_{\wh{Iso}(\fa)})=\wh{Iso}(\fa)/\wh{Fix}(\fa)=Iso(\fa)/Fix(\fa)\,.
\]
This proves theorem \ref{t:4}.
\qed
\subsection{Logarithmic foliations}\label{ss:33}

In this section we prove the results stated before concerning the isotropy group of logarithmic foliations.

\subsubsection{Proof of proposition \ref{p:4}}\label{ss:331}

Let $\om\in\Om^1(\C^n,0)$ be integrable and having an integrating factor $f\in\O_n$, so that,
$\Om:=f^{-1}.\,\om$ is closed. 
We will assume that, $\fa_\om$ has no meromorphic first integral and we want to prove that for any
$\Phi\in \wh{Iso}(\fa_\om)$ then $\Phi^*(\Om)=\d.\,\Om$, where $\d$ is a root of unity.

First of all, recall that
$\Phi^*(\om)=u.\,\om$ where $u\in\wh\O_n^*$. We assert that $\Phi^*(f)=v.\,f$, where $v\in\wh\O_n^*$

In fact, since $f^{-1}.\,\om$ is closed we have $\frac{df}{f}\wedge\om=d\om$. Applying $\Phi^*$ in both members of this relation we get
\[
\frac{d(f\circ\Phi)}{f\circ\Phi}\wedge\Phi^*(\om)=d(\Phi^*(\om))\,\,\iff\,\,\left(\frac{d(f\circ\Phi)}{f\circ\Phi}-\frac{du}{u}\right)\wedge\om=d\om
\]

If we set $g=\frac{f\circ\Phi}{u}$ then the above relation becomes
\[
\frac{dg}{g}\wedge\om=d\om\,\,\implies\,\,\left(\frac{dg}{g}-\frac{df}{f}\right)\wedge\om=0\,\,\implies\,\,d(g/f)\wedge\om=0\,.
\]
Therefore, $g/f$ is a meromorphic first integral of $\fa_\Om$ and so
$f/g=c$, where $c\in\C^*$. Hence,
\[
\Phi^*(\Om)=\Phi^*\left(\frac{\om}{f}\right)=C.\,\Om\,,
\]
where $C=1/c$.

It remains to prove that $C$ is a root of unity.
Let $f=\Pi_{j=1}^rf_j^{k_j}$, $k_j\ge1$, be the decomposition of $f$ into irreducible factors, so that
\[
\Om=\sum_{j=1}^r\la_j\frac{df_j}{f_j}+d\left(\frac{H}{f_1^{k_1-1}...f_r^{k_r-1}}\right)\,,
\]
where $\la_1,...,\la_r\in\C$. Note that:
\begin{itemize}
\item[i.] $\la_j\ne0$ for some $j$, because otherwise $\frac{H}{f_1^{k_1-1}...f_r^{k_r-1}}$ would be a meromorphic first integral. Without lost of generality, we will assume that $j=1$.
\item[ii.] $\Phi^*$ permutes the factors of $f$: there exists a permutation $\si\in S_r$ such that $f_j\circ\Phi=u_j.\,f_{\si(j)}$, where $u_j\in \wh\O^*_n$.
\end{itemize}
In particular, we have
\[
C.\,\Om=\Phi^*(\Om)=\sum_{j=1}^r\la_j\left(\frac{df_{\si(j)}}{f_{\si(j)}}+\frac{du_j}{u_j}\right)+d\Phi^*\left(\frac{H}{f_1^{k_1-1}...f_r^{k_r-1}}\right)
\]
Comparing the residues, we get
\[
\la_{\si(i)}=C.\,\la_i\,,\,1\le i\le r\,.
\]
Let $m\in\{1,...,r\}$ be such that $\si^m(1)=1$. From the above relation we get
\[
\la_1=\la_{\si^m(1)}=C.\,\la_{\si^{m-1}(1)}=C^2.\,\la_{\si^{m-2}(1)}=...=C^m.\,\la_1\,\implies\,C^m=1\,.\qed
\]

\begin{rem}\label{r:32}
\rm It follows from the proof of proposition \ref{p:4} that if $\Phi^*(f_j)=u_j.\,f_j$, $\forall\,1\le j\le r$, then $\d=1$: $\Phi^*(\Om)=\Om$.
\end{rem}

\subsubsection{Proof of theorem \ref{t:5}}\label{ss:332}
Let $\Om=\sum_{j=1}^r\la_j\frac{df_j}{f_j}$ satisfying the hypothesis of theorem \ref{t:5} and let
$\wt\Phi\in \wh{Iso}(\fa_\Om)$ be such that $D\wt\Phi(0)=\rho.\,I$, where either $\rho$ is not a root of unity, or the class of $\wt\Phi$ in $\wh{Iso}(\fa_\Om)/\wh{Fix}(\fa_\Om)$ has infinite order.
Since $\wt\Phi$ permutes the hypersurfaces $(f_j=0)$,
$1\le j\le r$, there exists $N\le r!$ such that $\wt\Phi^N(f_j=0)=(f_j=0)$, $1\le j\le N$. 
Note that $\Phi:=\wt\Phi^N$ fixes all hypersurfaces $(f_j=0)$ and $D\Phi(0)=\rho^N.\,I$.
\begin{lemma}\label{l:32}
Let $\phi\in \wh{Iso}(\fa_\Om)$ be such that $D\phi(0)=I$.
Then $\phi\in \wh{Fix}(\fa_\Om)$.
\end{lemma}

{\it Proof.}
We will consider the more general case where
$D\phi(0)=\be.\,I$.

Let $\Pi\colon(\wt\C^n,E)\to(\C^n,0)$ be a blow-up at $0\in\C^n$, where
$E\simeq\p^{n-1}$ is the exceptional divisor.
Denote by $\wt\fa_\Om$ the foliation $\Pi^*(\fa_\Om)$.

In the chart $U_1=\{(x,t)\in\C\times\C^{n-1}\,|\,x\in\C\,,\,t=(t_2,...,t_n)\in\C^{n-1}\}\sub \wt\C^n$ we have
\[
\Pi(x,t)=(x,x\,t_2,...,x\,t_n)=(x,x\,t)\in\C^n\,,
\]
and $E\cap U_1=(x=0)$.

Moreover, since the order of $f_j$ at $0\in\C^n$ is $k_j$, we can write
\[
f_j\circ\Pi(x,t)=x^{k_j}.\,\wt{f}_j(x,t)\,,
\]
where $\wt{f}_j$ is the strict transform of $f_j$ in this chart.
Note that:
\begin{itemize}
\item[iii.] The foliation $\wt\fa_\Om$ is defined in this chart by
\begin{equation}\label{eq:15}
\wt\Om:=\Pi^*(\Om)=\a.\,\frac{dx}{x}+\sum_{j=1}^r\la_j\frac{d\wt f_j}{\wt f_j}\,,
\end{equation}
where $\a=\sum_{j=1}^rk_j.\,\la_j\ne0$.
\end{itemize}
 In particular, $\wt\Om$ is non-dicritical and
$L_E:=E\setminus \bigcup_j(\wt{f}_j=0)$ is a leaf of $\wt\fa_\Om$.

Let $p=(0,t_o)\in L_E$ and fix. We assert that there exists a local chart
$(U,(\wt x,t)\in(\C,0)\times(\C^{n-1},0)$ such that
\[
\wt\Om|_U=\a\,\frac{d\wt x}{\wt x}\,.
\] 

In fact, since $\wt{f}_j(0,t_o)\ne0$, $1\le j\le r$, the form $\te:=\sum_j\la_j\frac{d\wt f_j}{\wt f_j}$ is exact in a neighborhood $U$ of $(0,t_o)$. In particular, there exists $h\in\O(U)$ such that $\te|_U=dh$. Therefore, if $\wt{x}=x.\,e^h$ then $\wt\Om|_U=\a\,\frac{d\wt x}{\wt x}$.

Now, there exists $\wt\phi\in \wh{Diff}(\wt\C^n,E)$ such that $\Pi\circ\wt\phi=\phi\circ\Pi$. The reader can check that in the chart $(U,(\wt x,t))$ the transformation $\wt\phi$ is of the form
\[
\wt\phi(\wt x,t)=(g_1(\wt x,t),g_2(\wt x,t),...,g_n(\wt x,t))=(g_1(\wt x,t),t+\wt x.\,H(\wt x,t))\,
\]
where $H\in \O(U)$. In particular, $\wt\phi|_E=id_E$, the identity map of $E$.

On the other hand, by proposition \ref{p:4} we have $\phi^*(\Om)=\d\,\Om$, where $\d$ is a root of unity. From $\Pi\circ\wt \phi=\phi\circ\Pi$ we get
\[
\wt \phi^*(\wt\Om)=\wt \phi^*\circ\Pi^*(\Om)=(\Pi\circ\wt\phi)^*\circ(\Om)=\Pi^*\circ\phi^*(\Om)=\Pi^*\d\,\Om=\d\,\wt\Om\,\implies
\]
\[
\wt\phi^*\left(\a.\frac{d\wt x}{\wt x}\right)=\d\,\a.\,\frac{d\wt x}{\wt x}\,\implies\,\d=1\,\,\text{and}\,\,
\frac{d g_1}{g_1}=\frac{d\wt x}{\wt x}\,\implies\,g_1(\wt x,t)=\be.\,\wt x\,.
\]
In the chart $(\wt x,t)$ the leaves of $\wt\fa_\Om$ are of the levels $\wt x=$ constant.
If $\be=1$ then $\wt\phi(\wt x,t)=(\wt x,t+x.\,H(\wt x,t))$, so that $\wt\phi$ preserves the leaves of $\wt\fa_\Om$, which implies that $\phi\in \wh{Fix}(\fa_\Om)$.
\qed

\begin{cor}\label{c:31}
Let $\Phi\in\wh{Iso}(\fa_\Om)$ such that $D\Phi(0)=\rho.\,I$, where $\rho\ne1$. Then the class of $\Phi$ in $\wh{Iso}(\fa_\Om)/\wh{Fix}(\fa_\Om)$ has a formally linearizable representative $\wt\Phi$. 
\end{cor}

{\it Proof.} Let $\Phi\in\wh{Iso}(\fa_\Om)$ be such that $D\Phi(0)=\rho.\,I$. Let $\Phi=\Phi_S\circ\Phi_U$ be the decomposition of $\Phi$ as in proposition \ref{p:22}. Recall that $\Phi_S,\Phi_U\in\wh{Iso}(\fa_\Om)$, $\Phi_S$ is formally linearizable and $\Phi_U$ is unitary. Since
$\rho.\,I=D\Phi_S(0)\circ D\Phi_U(0)$ we have $D\Phi_S(0)=\rho.\,I$ and $D\Phi_U(0)=I$. Therefore, by lemma \ref{l:32} $\Phi_U\in \wh{Fix}(\fa_\Om)$ and the classes of $\Phi_S$ and $\Phi$ in $\wh{Iso}(\fa_\Om)/\wh{Fix}(\fa_\Om)$ are equal. This proves corollary \ref{c:31}.
\qed

\begin{cor}\label{c:32}
Let $\Phi\in \wh{Iso}(\fa_\Om)$ with $D\Phi(0)=\rho.\,I$. If the class of $\Phi$ in
$\wh{Iso}(\fa_\Om)/\wh{Fix}(\fa_\Om)$ has infinite order then $\rho$ is not a root of unity.
Moreover:
\begin{itemize}
\item[(a).] $\Phi$ is formally linearizable: there exists $g\in \wh{Diff}(\C^n,0)$ such that
$g^{-1}\circ \Phi\circ g=D\Phi(0)$.
\item[(b).] $g^*(\Om)=\sum_{j=1}^r\la_j\frac{dh_j}{h_j}$. In particular $g^*(\Om)$ is homogeneous.
\end{itemize}
\end{cor}

{\it Proof.} 
It follows from lemma \ref{l:32} that $\rho\ne1$. Therefore by corollary \ref{c:31} there is a formally linearizable representative of $\Phi$ in $\wh{Iso}(\fa_\Om)/\wh{Fix}(\fa_\Om)$, say $\wt\Phi$. If $\rho^k=1$, for some $k\in\N$ then $\wt\Phi^k=I$, so that its class has finite order. Therefore, $\rho$ is not a root of unity.

In particular, $D\Phi(0)$ is non-resonant and $\Phi$ is formally linearizable by Poincar\'e's linearization theorem: there exists an unique $g\in \wh{Diff}(\C^n,0)$ such that $Dg(0)=I$ and 
$g^{-1}\circ \Phi\circ g=\rho.\,I:=T\in GL(n,\C)$.

Let us prove (b). Set $\Om^*:=g^*(\Om)$, so that $T^*\Om^*=\Om^*$ by proposition \ref{p:4}.
On the other hand, the Zariski closure of the group $H=\{T^n\,|\,n\in\Z\}$ is the group
$\C^*.\,I\sub GL(n,\C)$. This implies that, if $T_\tau=e^\tau.\,I=exp(\tau\,R)$ then
\[
exp(\tau\,R)^*\,\Om^*=\Om^*\,,\,\forall\tau\in\C^*\,\implies\,L_R\,\Om^*=\frac{d}{d\tau}(exp(\tau\,R)^*\,\Om^*)|_{\tau=0}=0\,.
\]
Let us finish the proof that $\Om^*=\sum_{j=1}^r\la_j\frac{dh_j}{h_j}$.
Set $f_j^*:=f_j\circ g\in\wh\O_n$, $1\le j\le r$, so that
$\Om^*=\sum_{j=1}^r\la_j\frac{df_j^*}{f_j^*}$. Note that $j_0^{k_j}f_j^*=h_j$, because $Dg(0)=I$. Now,
\[
L_R\Om^*=i_R\,d\Om^*+d\,i_R\Om^*=d\,i_R\Om^*\,\implies\,d\,i_R\Om^*=0\,\implies\,i_R\Om^*=c\,\iff
\]
\[
\sum_{j=1}^r\la_j\,f_1^*...R(f_j^*)...f_r^*=c\,f_1^*...f_r^*\,\implies\,R(f_j^*)=v_j.\,f_j^*\,,\,v_j\in\wh\O_n^*\,,
\]
because $f_j^*$ is irreducible, $1\le j\le r$. As the reader can check, the last relation implies that $h_j$ divides $f_j^*$: $f_j^*=\wt{v}_j.\,h_j$. Therefore,
\[
\Om^*=\sum_{j=1}^r\la_j\frac{dh_j}{h_j}+\Te\,,\,\Te=\sum_j\la_j\frac{d\wt{v}_j}{\wt{v}_j}\in\wh\Om^1(C^n,0)\,.
\]
Since $h_j$ is homogeneous of degree $k_j$ we have $R(h_j)=k_j\,h_j$ by Euler's identity, so that
$L_R\left(\frac{dh_j}{h_j}\right)=0$. This implies that $L_R\Te=0$ and so $\Te=0$.
\qed
\vskip.1in
If $\Phi$ is holomorphic and $|\rho|\ne1$ or if $|\rho|=1$ and $\rho$ satisfies a small denominator condition then
$g\in Diff(\C^n,0)$ (cf. \cite{ad}) and we are done.

In the general case however, when $\Phi$ is only formal, the idea is to prove that there exists a holomorphic vector field $\wt{R}\in \X(\C^n,0)$ with $D\wt{R}(0)=R$, the radial vector field,
with $i_{\wt R}\Om=\a=\sum_jk_j.\,\la_j$.
In the proof of this fact, we will use Artin's approximation theorem.

Let us finish the proof of theorem \ref{t:5}, assuming the existence of a vector field $\wt R$ as above.
If $i_{\wt R}\Om=\a$ then $L_{\wt R}\Om=i_{\wt R}d\Om+d\,i_{\wt R}\Om=0$,
because $\Om$ is closed and $i_{\wt R}\Om$ is a constant.

On the other hand, since $D\wt R(0)=R$, by Poincar\'e's linearization theorem (cf. \cite{ad}), the vector field $\wt R$ is holomorphically linearizable: there exists $\phi\in Diff(\C^n,0)$ such that
$D\phi(0)=I$ and $\phi^*(\wt R)=R$.  
We assert that $\phi^*(\Om)$ is homogenous:
\[
\phi^*(\Om)=\sum_{j=1}^r\la_j\frac{dh_j}{h_j}\,.
\]
In fact, set $\wt{f}_{j}=\phi^*(f_j)$ and $\wt\Om=\phi^*\Om$. Then
\[
\wt\Om=\sum_{j=1}^r\la_j\frac{d\wt{f}_j}{\wt{f}_j}\,\,\text{and}\,\,i_R\wt\Om=\a\,\implies
\]
$R(\wt{f}_j)=u_j\,\wt{f}_j$, $u_j\in\O_n$ and $u_j(0)=k_j$. Writing the Taylor series of $\wt{f}_j$ and of $u_j$, by an induction argument, we obtain that there exists an unity $v_j\in\O_n^*$ such that $\wt{f}_j=v_j.\,h_j$.
In particular, we get
\[
\Om_1=\sum_{j=1}^r\la_j\frac{dh_j}{h_j}+\Te\,,\,\text{where}\,\Te=\sum_j\la_j\frac{dv_j}{v_j}\in\Om^1(\C^n,0)\,.
\]
Since $h_j$ is homogeneous, $1\le j\le r$, we have
\[
L_R\left(\sum_j\la_j\frac{dh_j}{h_j}\right)=0\,\implies\, L_R\Te=0\,\implies\,\Te=0\,,
\]
because $\Te$ is holomorphic. Therefore, $\wt\Om$ is homogeneous.
\vskip.1in
It remains to prove the existence of $\wt R$ with $i_{\wt R}\Om=\a$.
We have proved that there exists $g\in\wh{Diff}(\C^n,0)$ such that
\[
g^*(\Om)=\sum_{j=1}^r\la_j\frac{dh_j}{h_j}:=\Om^*\,.
\]
In particular, if $R$ is the radial vector field on $\C^n$ then 
\[
i_R\Om^*=\sum_{j=1}^r\la_j\frac{R(h_j)}{h_j}=\a\,.
\]

Let $\wh R=g_* R\in \wh\X(\C^n,0)$ and note that $i_{\wh R}\Om=\a$, because $i_Rg^*\Om=\a$.
In fact, let $\wh R=\sum_{j=1}^n\var_j\frac{\pa}{\pa z_j}$, where $\var_1,...,\var_n\in \wh\O_n$. Writing explicitly the relation $i_{\wh R}\Om=\a$ we get
\[
\sum_{j=1}^r\la_j\frac{\wh R(f_j)}{f_j}=\a\,\implies\,\sum_{j=1}^rf_1...\wh{R}(f_j)...f_r=\a.\,f_1...f_r\,\implies
\]
\[
\implies\,\underset{1\le i\le n}{\sum_{1\le j\le r}}\var_j.\,f_1...\frac{\pa f_j}{\pa z_i}...f_r=
\a.\,f_1...f_r\,\,.
\]
In particular, $\var=(\var_1,...,\var_n)\in\wh\O_n^{\,n}$ is a formal solution of the analytic equation $F(z,w)=0$, where
\[
F(z,w)=\underset{1\le i\le n}{\sum_{1\le j\le r}}w_j.\,f_1(z)...\frac{\pa f_j(z)}{\pa z_i}...f_r(z)-
\a.\,f_1(z)...f_r(z)\,.
\]
It follows from Artin's approximation theorem that $F(z,w)=0$ has a convergent solution $w=\phi(z)=
(\phi_1(z),...,\phi_n(z))$ such that $J^1_0(\phi)=J^1_0(\var)$. Since $j^1_0(\wh R)=R$, the radial vector field, we can conclude that the vector field $\wt R=\sum_{j=1}^n\phi_j\,\frac{\pa}{\pa z_j}$ satisfies $i_{\wt R}\Om=\a$ and $D\wt R(0)=R$, as desired.
This finishes the proof of theorem \ref{t:5}.
\qed

\subsubsection{Proof of corollaries \ref{c:3}, \ref{c:4} and \ref{c:5}}
$\,$
\vskip.1in
{\it Proof of corollary \ref{c:3}.} 
Let us sketch how theorem \ref{t:5} implies corollary \ref{c:3}. Let $\wt\Phi\in \wh{Iso}(\fa_\Om)$ whose class in $\wh{Iso}(\fa_\Om)/\wh{Fix}(\fa_\Om)$ has infinite order. Since $\fa_\Om$ has no formal meromorphic first integral, then $\wt\Phi$ permutes the set of hypersurfaces $(f_j=0)$, $1\le j\le r$. In particular, there exists $N\in\N$ such that $f_j\circ\wt\Phi^N=u_j.\,f_j$, where $u_j\in\wh\O_n$. If $T=D\wt\Phi^N(0)$ then $h_j\circ T=u_j(0).\,h_j$, $\forall$ $1\le j\le r$. Since $(h_1,...,h_r)$ is rigid we must have $D\wt\Phi^N(0)=\rho.\,I$. In particular $\Phi:=\wt\Phi^N$ satisfies hypothesis (a) of theorem \ref{t:5}, as wished. This finishes the proof of corollary \ref{c:3}.
\qed
\vskip.1in
{\it Proof of corollary \ref{c:4}}. 
First of all we would like to observe that there exists $N\in\N$ such that for any $\Phi\in \wh{Iso}(\fa_\Om)$ then $D\Phi^N(0)=\rho\,I$, where $\rho$ is a root of unity.

In fact, there exists $N\in \N$ with $N\le r!$ and such that $f_j\circ \Phi^N=u_j.\,f_j$, $1\le j\le r$.
Since the set $(h_1,...,h_r)$ is rigid, we have $D\Phi^N(0)=\rho.\,I$.
If $\rho$ was not a root of unity then by theorem \ref{t:5} there would exist
$\psi\in \wh{Diff}(\C^n,0)$ such that $f_1\circ \psi=u_1.\,h_1$. However, this is impossible, because $f_1$ is irreducible and $h_1$ is not.

\begin{claim}\label{cl:31}
\rm Let $\Si=\{k\in\N\,|\,$ there exists a primitive $k^{th}$-root of unity $\rho$ such that $D\Phi(0)=\rho.\,I$, for some $\Phi\in\wh{Iso}(\fa_\Om)$ such that $f_j\circ\Phi=u_j.\,f_j$, $\forall j\,\}$\,. We claim that $\Si$ is bounded.
\end{claim}

{\it Proof.} We will use that $f_1$ is irreducible and $h_1$ is not.
First of all, we assert that for any $\rho\in X$ there exists $\Phi\in\wh{Iso}(\fa_\Om)$ such that
$D\Phi(0)=\rho.\,I$ and $\Phi$ is formally linearizable.

Fix $\ell\in \Si$, so that there is $\Phi\in\wh{Iso}(\C^n,0)$ with $D\Phi(0)=\rho\,I$, where $\rho$ is a $\ell^{th}$-root of unity and $f_j\circ\Phi=u_j.\,f_j$, $1\le j\le r$.
Let $\Phi=\Phi_S\circ\Phi_U$, where $\Phi_S$ is semi-simple, $\Phi_U$ unitary and $\Phi_S\circ\Phi_U=\Phi_U\circ\Phi_S$.
As we have seen in remark \ref{p:22}, $\Phi_S,\Phi_U\in \wh{Iso}(\fa_\Om)$.
On the other hand, $D\Phi_U(0)=I$ and $D\Phi_S(0)=\rho.\,I$. Since $\Phi_S$ is linearizable, this proves the assertion.

Let $\ell\in \Si$ and $\Phi\in \wh{Fix}(\fa_\Om)$ be formally linearizable with $D\Phi(0)=\rho.\,I$.
Since $\Phi$ is linearizable there exists $\psi\in \wh{Diff}(\C^n,0)$ such that $\psi^{-1}\circ \Phi\circ\psi=\rho.\,I$.
Let $f_1^*=f_1\circ\psi$. We assert that, if $\ell>k_1$ then $h_1$ divides the $(\ell+k_1-1)^{th}$-jet of $f_1^*$.

In fact, let $f_1^*=\sum_{i\ge k_1}g_i$ be the Taylor series of $f_1^*$, where $g_i$ is homogeneous of degree $i$. Note that $g_{k_1}=h_1$, because $D\psi(0)=I$.
Since $f_1\circ \Phi=u_1.\,f_1$, where $u_1\in \wh\O_n^*$, we have
\begin{equation}\label{eq:18}
f_1^*\circ(\rho.\,I)=u_1^*.\,f_1^*\,,
\end{equation}
where $u_1^*=u_1\circ\psi$. Let $u_1^*=\sum_{i\ge0}w_i$ be the Taylor series of $u_1^*$. Relation 
(\ref{eq:18}) can be written as
\[
\sum_{i\ge k_1}\rho^i\,g_i=\sum_{i\ge k_1}\sum_{s=k_1}^i  w_{i-s}.\,g_s\,\implies\,w_0=\rho^{k_1}\,\,\text{and}\,\,
(\rho^i-\rho^{k_1})\,g_i=\sum_{s=k_1}^{i-1}w_{i-s}\,g_s\,,\forall i>k_1
\]
Since $\rho$ is a primitive $\ell^{th}$ root of unity, $\rho^i-\rho^{k_1}\ne0$ if $k_1+1\le i\le \ell+k_1-1$. Therefore, the above relation implies that $h_1|g_i$ if $k_1\le i\le \ell+k_1-1$, as the reader can check. Hence, $h_1$ divides the $(\ell+k_1-1)^{th}$ jet of $f_1^*$, as asserted.
 
Now, suppose by contradiction that $\Si$ was unbounded. In this case, since $h_1$ is reducible, the above argument implies that for any $k\in \N$ there exists $\psi\in \wh{Diff}(\C^n,0)$ such that
$J^k_0(f_1\circ \psi)$ is reducible. On the other hand, we have
\[
J^k_0(f_1\circ \psi)=J^k_0(J^k\,f_1\circ J^k_0\,\psi)\,\implies\,J^k_0f_1\,\text{is reducible
$mod\,m_n^{k+1}$}\,,
\]
in the sense that there exist $f,g\in \O_n$ such that $f(0)=g(0)=0$ and
$J_0^kf_1=J_0^k(f.\,g)$. Therefore, for any $k\in\N$ then $J_0^k(f_1)$ is reducible $mod\,m_n^{k+1}$, and this implies that $f_1$ is reducible.
This proves the claim,
\qed
\vskip.1in
Now, let $\De\colon \wh{Iso}(\fa_\Om)\to GL(n,\C)$ be the homomorphism $\De(\Phi)=D\Phi(0)$.
Note that lemma \ref{l:32} implies that $ker(\De)\sub \wh{Fix}(\fa_\Om)$.
We assert that, in fact $ker(\De)=\wh{Fix}(\fa_\Om)$. 

In fact, if $\Phi\in \wh{Fix}(\fa_\Om)$ then $f_j\,|\,f_j\circ\Phi$, $1\le j\le r$, because $(f_j=0)$ is a leaf of $\fa_\Om$. This implies that $D\Phi(0)\in \I(h_1,...,h_r)$ (see definition \ref{d:3}). Since $(h_1,...,h_r)$ is rigid, we $\De(\Phi)=I$. Hence, $ker(\De)=\{I\}$.

This implies that $\wh{Iso}(\fa_\Om)/\wh{Fix}(\fa_\Om)\simeq Im(\De)$.
On the other hand, since $X$ is bounded, the subgroup $G:=\C^*\,I\cap Im(\De)$, of $Im(\De)$, is finite. This subgroup is contained in the center of $Im(\De)$,  and so is a normal subgrup.

In order to conclude the proof of corollary \ref{c:4} is sufficient to prove that $Im(\De)/G$ is finite. However, this group can be identified with a subgroup of the group $S_r$ of permutations of $\{1,...,r\}$.
 
In fact, as we have seen before, if $\Phi\in \wh{Iso}(\fa_\Om)$ then there exists a permutation $\si\in S_r$ such that
$f_j\circ\Phi=u_j.\,f_{\si(j)}$, $\forall j$. Moreover, $\Phi\in G$ if, and only if, $\si=id$, which proves the assertion.
This finishes the proof of corollary \ref{c:4}.
\qed


\bibliographystyle{amsalpha}

\vskip.3in

{\sc D. Cerveau}

{\sl Universit\'e de Rennes, CNRS,}

{\sl IRMAR-UMR 6625,}

{\sl F-35000 Rennes, France}

{\tt E-Mail: dominique.cerveau@univ-rennes1.fr}

\vskip.3in

{\sc A. Lins Neto}

{\sl Instituto de Matem\'atica Pura e Aplicada}

{\sl Estrada Dona Castorina, 110}

{\sl Horto, Rio de Janeiro, Brasil}

{\tt E-Mail: alcides@impa.br}

\end{document}